\documentclass[11pt,a4paper]{article}

\usepackage{latexsym}
\usepackage{graphicx}
\usepackage{times}
\usepackage{psfrag}
\usepackage{bbm}
\usepackage{bm}
\usepackage{setspace}
\usepackage{url}
\usepackage{comment}
\usepackage{color}

\usepackage[round]{natbib}

\renewcommand{\vec}{\bm}
\newcommand{\dotp}[1]{\left\langle #1\right\rangle}
\newcommand{\norm}[1]{\| #1 \|}
\newcommand{\floor}[1]{\left\lfloor #1\right\rfloor}
\newcommand{\ceil}[1]{\left\lceil #1\right\rceil}

\newcommand{\RR}{\mathbbm{R}}
\newcommand{\ZZ}{\mathbbm{Z}}

\newcommand{\var}{\mathrm{var}}

\def\Chi{%
    \mbox{
    {\kern-.20em\setbox0=\hbox{X}%
    \vbox to 1.0\ht0{\hbox{$\chi$}\vss}}%
    \kern-.00em} }
\def\mathlatex{%
    \mbox{
    L\kern-.36em {\setbox0=\hbox{T}%
    \vbox to \ht0{\hbox{\the\scriptfont0 A}\vss}}%
    \kern-.15em \TeX} }

\newcommand{\barbar}[1]{\overline{\overline{#1}}}

\newcommand{\proof}{\textbf{Proof. }}
\newcommand{\qed}{\hbox{}\hfill $\Box$ }
\newtheorem{lemma}{Lemma}
\newtheorem{thm}{Theorem}

\newtheorem{definition}{Definition}
\newtheorem{corollary}{Corollary}
\newtheorem{prop}{Proposition}
\newtheorem{remark}{Remark}
\newenvironment{shortversion}{}{}
\newenvironment{longversion}{}{}
\newenvironment{revision}{}{}

\excludecomment{shortversion}

\begin{shortversion}
\title{Non-equispaced B-spline wavelets}
\end{shortversion}
\begin{longversion}
\title{Non-equispaced B-spline wavelets\\
\fbox{\normalsize{\textbf{Version with detailed proofs}}}}
\end{longversion}
\author{Maarten Jansen \\
Universit\'e libre de Bruxelles\\ 
Departments of Mathematics and Computer Science}

\begin{document}

\maketitle

\begin{abstract}
This paper has three main contributions. The first is the construction of
wavelet transforms from B-spline scaling functions defined on a grid of
non-equispaced knots. The new construction extends the equispaced, 
biorthogonal, compactly supported Cohen-Daubechies-Feauveau wavelets. The new
construction is based on the factorisation of wavelet transforms into lifting
steps. The second and third contributions are new insights on how to use these
and other wavelets in statistical applications. The second contribution is
related to the bias of a wavelet representation. It is investigated how the
fine scaling coefficients should be derived from the observations. In the
context of equispaced data, it is common practice to simply take the
observations as fine scale coefficients. It is argued in this paper that this
is not acceptable for non-interpolating wavelets on non-equidistant data.
Finally, the third contribution is the study of the variance in a
non-orthogonal wavelet transform in a new framework, replacing the numerical
condition as a measure for non-orthogonality. By controlling the variances of 
the reconstruction from the wavelet coefficients, the new framework allows us
to design wavelet transforms on irregular point sets with a focus on their use
for smoothing or other applications in statistics.
\end{abstract}

\subsubsection*{Keywords}
B-spline; wavelet; Cohen-Daubechies-Feauveau; penalised splines;
non-equispaced; non-equidistant; lifting

\subsubsection*{AMS classification}
42C40; 65T60; 65D07; 41A15; 65D10


\section{Introduction}
\label{sec:intro}

Ever since the early days of wavelet research, spline wavelets have enjoyed
special attention in the community. Spline wavelets combine the benefits from a
sparse multiscale approach using wavelets and the well known properties of
splines, including the closed form expressions, the numerous recursion
relations, the polynomial based regularity of the basis functions
\citep{unser97:tengoodreasons}.

\begin{revision}
Splines \citep{deBoor01:splines}, formally defined by the recursion as in 
(\ref{eq:Bsplinerecursion}) below, are piecewise polynomials of a certain
degree, with continuity constraints in the knots that connect the polynomial
pieces. The position of these knots and the degrees of the polynomial pieces
are key parameters in the numerous methods in computer aided geometric design
and computer graphics that are based on splines.
Splines are also a popular
tool in numerical analysis, for instance in interpolation. Compared to full
polynomial interpolation, spline interpolation is far less sensitive to
numerical instabilities that lead to oscillations. The good numerical condition
is linked to the fact that any spline function can be decomposed into a basis 
of compactly supported piecewise polynomials, so-called B-splines. In
statistics, the coefficients of the B-spline decomposition can be estimated in
a nonparametric regression context. The estimator typically minimizes the
residual sum of squares, penalized by the roughness of the regression curve
\citep{GreenSilverman94,eubank99:splinesmoothing}. 
The spline wavelet smoothing, as discussed in Section
\ref{sec:estBsplinewavelets} of this paper, can be considered as an extension
of these smoothing splines towards sparsity oriented penalties and
corresponding nonlinear smoothing based on thresholding. While smoothing
splines have their knots on the locations of the observations,
P-splines \citep{ruppert03:semiparametric} allow a flexible choice of knots.
An important advantage of any spline, whether it be an interpolating, smoothing
or P-spline, is that it is known by an explicit expression. The main merit
of a closed-form expression is that all information about the smoothness
of the function, typically expressed by the Lipschitz regularity, is readily
available for use in smoothing algorithms.
The spline wavelets on irregular knots constructed in this 
paper are splines, and thus share this benefit. This is in contrast to most
other wavelets, especially those on irregular point sets. The smoothness of
these wavelets depends on the limit of an infinitely iterated refinement or
subdivision scheme \citep{daubechies99:regir}, which can be hard to analyse, 
even for straightforward refinement schemes such as for Deslauries-Dubuc
interpolating wavelets 
\citep{deslauriers87:interpolation,deslauriers89:symmetric,donoho99:deslauries,sweldens96:athome}.
\end{revision}

The combination of splines and wavelets is, however, not trivial. One of the
problems is that splines do not provide a naturally orthogonal basis. On the
other hand, orthogonality is much appreciated in wavelet theory, because
wavelet transforms operate scale after scale. 
\begin{revision}
In the inverse transform, i.e., in the reconstruction of a function from its
wavelet coefficients, this scale to scale process amounts to the refinement
or subdivision, mentioned above. Although the smoothness of the reconstruction
in a spline basis does not depend on this refinement scheme, other properties
do. These properties include the numerical condition of the transform as well
as the bias and the variance in statistical estimation. The assumption of
orthogonality facilitates the analysis of these properties throughout the
subdivision scheme.
\end{revision}
Initial spline wavelet constructions were orthogonal
\citep{battle87:ondelettes,lemarie88:ondelettes}.
The price to pay for the orthogonality was that the basis functions did not
have a compact support, and related to this, that transformation matrices
were full, not sparse, matrices. 

The condition of orthogonality can be 
relaxed if the basis functions within one scale are allowed to be
non-orthogonal, while the wavelets at different scales are still kept 
orthogonal
\citep{chui92:splinewavelets,unser92:splinewavelets,unser93:familysplinewavelets}.
This construction leads to a semi-orthogonal spline basis. Since the
non-orthogonality occurs only within each scale, this has little impact on the
asymptotic analysis of the refinement process. Moreover, semi-orthogonal
wavelet bases can be constructed from non-orthogonal spline bases, such as
B-splines. The B-splines and the resulting wavelets have compact support. As a
consequence, the reconstruction of data from a decomposition in these bases
uses a sparse matrix. Sparse matrices and compact supports lead to faster
algorithms, but also contribute to better control over manipulations on the
wavelet coefficients. A manipulation of a coefficient, such as a thresholding
or shrinkage operation, has only a local effect.

Unfortunately, the forward transform of observations into the semi-orthogonal
wavelet basis still requires the application of a full, non-sparse, matrix. 
\begin{revision}
This is because in general the inverse of a sparse matrix is not sparse. Sparse
inverses are possible in a wide range of fast wavelet transforms, but the
combination of the semi-orthogonality and the spline properties cannot
be obtained within a sparse forward transform.
\end{revision}
As a consequence, every fine scale observation has at least some contribution
to each wavelet coefficient. It would be better if not only a wavelet
coefficient had local impact on the reconstruction of the observations but
also, at the same time, the coefficient got its value from a limited number
of observations. The latter, dual, form of compact support is possible in the
framework of bi-orthogonal spline wavelets. The construction by 
\citet*{coh-dau-fea92:biorthogonal}
of a multiresolution
analysis starting from B-splines has led to non-orthogonal wavelets with sparse
decomposition and reconstruction matrices.

The Cohen-Daubechies-Feauveau wavelets are defined on equispaced knots. This is
because the classical multiresolution theory starts from scaling bases whose
basis functions are all dilations and translations of a single father function.
This construction is not possible on irregular knots. 
B-splines, on the other hand, are easily defined on non-equispaced knots.
This paper extends the construction by Cohen, Daubechies and Feauveau towards
these non-equispaced B-splines. For the construction of wavelet transforms on
non-equispaced knots, sometimes termed second generation wavelets
\citep{daubechies99:irregular,sweldens97:lifting2nd}, this paper adopts
the lifting scheme \citep{sweldens96:liftingbior}. The lifting scheme provides
for every refinement operation in a wavelet transform a factorisation
into simple steps within the scale of that refinement. The key contribution of
this paper is to identify the lifting steps that are necessary in the
factorisation of a B-spline refinement on non-equidistant knots.

Section \ref{sec:resultsfromliterature} summarizes results from the literature
that are necessary for the main contribution in Section
\ref{subsec:designliftingsteps}.
First, Section \ref{subsec:multilevelgrids} defines the notion of multiscale
grids. Then, Section \ref{subsec:Bsplines1scale} gives a definition of 
B-splines together with their properties for further use.
Section \ref{subsec:Bspline2scales} proposes to use a refinement equation as a
definition for B-splines. In order to fill in coefficients in the equation, it
needs to be factored into elementary operations. The result is a
slight generalisation of the well known factorisation of equidistant wavelet
transforms \citep{daubechies98:factorlifting}.
Finally, Section \ref{subsec:multiscaletransform} further investigates the role
of the factorisation in wavelet transforms.
The main contribution in Section \ref{subsec:designliftingsteps} fills in the 
lifting steps that constitute a B-spline refinement.
Section \ref{subsec:Bsplinewavelets} gives an expression for all possible
wavelets that fit within the B-spline refinement scheme.
The first part of this paper, about the construction of non-equispaced B-spline
wavelets is concluded by the short Section \ref{subsec:RWT} on the
non-decimated B-spline wavelet transform.

Other work on lifting for spline or B-spline wavelets, such as
\citep{bertram04:liftingBsplines,chern99:CDFlifting,fahmy08:Bsplinelifting,li05:Bsplinewavelets,prestin05:splinelifting,xiang07:waveletfinelem}
is situated on equidistant knots or is focused on
specific cases, such as Powell-Sabin spline wavelets
\citep{vanraes04:powellsabine}.
B-spline wavelets on non-equispaced knots have been studied for specific
applications and with particular lifting schemes
\citep{pan09:nonuniformBsplinewavelets,lyche01:nonuniformsplinewavelets}.
It should be noted that B-splines on non-equispaced knots are often termed
non-uniform B-splines. This term is avoided in this paper, as in statistical
sense, a uniform set of knots could be interpreted as a set of random,
uniformly distributed knots, which are of course almost surely non-equidistant.

The second part of the paper consists of the Sections \ref{sec:apprxerr}
and \ref{sec:estBsplinewavelets}. 
It concentrates on the use of the B-spline wavelets from the first part in
statistics.
In statistical applications, B-spline wavelets are used, for instance, in a
soft threshold scheme for noise reduction. Given that a soft threshold comes
from an $\ell_1$ regularised least squares approximation of the input, this
application is an example of a penalised spline method.
The discussions in Sections \ref{sec:apprxerr} and \ref{sec:estBsplinewavelets}
are quite general, and therefore applicable to other non-equispaced wavelets as
well.

The discussion in Section \ref{sec:apprxerr} is related to the bias in a
B-spline wavelet smoothing. It investigates how to proceed from observations to
fine scaling coefficients. It is argued that in a situation with 
non-equispaced data, this step should be taken with care, in order to avoid to
commit what some authors call the ``wavelet crime''
\citep{strang96:filterbanks}.

While Section \ref{sec:apprxerr} deals with bias,
Section \ref{sec:estBsplinewavelets} is about the variance in a second
generation wavelet transform. Assuming that the observations are independent,
the variance of the transformed data is best understood if the transform is
orthogonal. As the construction by Cohen, Daubechies and Feauveau, extended in
this paper towards non-equidistant knots, has somehow less attention for
orthogonality, this may cause major problems with estimators suffering from
large variance effects
\citep{vanraes02:stabilizing,vanaerschot06:adaptivesplitting,jansen09:irregular}.
Although the large variance is due to the transform
being non-orthogonal, the classical numerical condition number is not a
satisfactory quantification of the statistical problem. Therefore this paper
proposes a multiscale variance propagation number, based on the singular values
of the linear projection onto the coarse scale B-spline basis. From there, an
alternative B-spline wavelet transform is developed. The alternative is closely
related to the Cohen-Daubechies-Feauveau construction, but it keeps the
variance propagation under control.


\section{B-splines and multiresolution}
\label{sec:resultsfromliterature}


This section reviews definitions and well established results on B-splines,
multiscale representations and the lifting scheme. 

\subsection{Multilevel grids}
\label{subsec:multilevelgrids}

Let $K_n = \{x_k|k = 0,\ldots,n-1\}$ be a set of knots on which we will define
B-splines. B-splines of order $\widetilde{p}$ are basis functions spanning
all piecewise polynomials of degree $\widetilde{p}-1$ with continuous
$\widetilde{p}-2$ derivatives in the knots. In this paper, the B-spline basis
will be constructed through a process known as refinement or subdivision. For
this process to work, we first have to define coarse scale versions of the grid
of knots. We thus identify the input set of knots as the fine scale grid,
formalised as $x_{J,k} = x_k$. The index $J$ refers to the highest or finest
scale. From there, we define grids $\{x_{j,k}|k = 0,\ldots,n_j-1\}$
at coarser scales $j$, where $j=L,L+1,\ldots,J-1$, $L$ being the lowest or
coarsest scale. Obviously $n_j < n$ stands for the number of points at scale
$j$. Denoting $\Delta_{j,k} = x_{j,k+1}-x_{j,k}$, and $\Delta_j =
\sup_{k=0,\ldots,n_j-2} \Delta_{j,k}$, we call the grid at scale $j$
\emph{regular} if $\Delta_{j,k}$ does not depend on $k$, i.e., $\Delta_{j,k} =
\Delta_j$. The focus in this paper lies on irregular grids. 

\begin{definition}(multilevel grid)
The sequence of grids constitutes a \emph{multilevel grid}
if the following conditions are met:
\begin{enumerate}
\item
The sequence $n_j$ is strictly increasing.
\item
There exist constants $R \in \RR$ and $\beta>0$ so that maximum gap at scale
$j$ is bounded as follows
\begin{equation}
\Delta_j \leq Rn_j^{-\beta}.
\label{eq:maxgapdecay}
\end{equation}
\end{enumerate}
\label{def:multilevelgrid}
\end{definition}
Condition (\ref{eq:maxgapdecay}) is a slightly stricter version of the
definition adopted in \citep{daubechies01:commutation}, where in the context of
binary or dyadic refinement, i.e., $n_j = 2^j$, it is imposed that
$\sum_{j=L}^\infty \Delta_j < \infty$. The condition can be understood by
considering a sequence of functions $x_j:[0,1]\to \RR:u\mapsto x_j(u)$ 
for which $x_{j,k} = x_j(k/(n_j-1))$. The divided differences of $x_j(u)$ in
the knots $k/(n_j-1)$ are then $\Delta_{j,k}(n_j-1)$. 
The divided differences must not or at most very slowly converge to a locally 
infinite derivative, in order not to leave any coarse scale gaps in a grid at
fine scale.

A multilevel grid is \emph{nested} if
$x_{j+1,k} \in \{x_{j,k}| k = 0,\ldots,n_j-1\}.$
In particular, the multilevel grid is two-nested if at each level, the grid is
a binary refinement of the previous, coarser level, that is, if $x_{j+1,2k} =
x_{j,k}$. This paper works with two-nested multilevel grids only.

\subsection{B-splines at a fixed scale}
\label{subsec:Bsplines1scale}

Throughout this paper, $\varphi_{j,k}^{\left[\widetilde{p}\right]}(x)$ will
stand for the B-spline of order $\widetilde{p}$ defined on the knots $x_{j,i}$.
There exist several recursion expressions for the construction of B-splines.
This paper will use the following formula
\citep[page 497]{qu92:subdivisionBsplines,daubechies01:commutation}
as definition.
\begin{definition}(B-splines)
The B-splines of order 1 defined on the knots $x_{j,i}$ are the characteristic
functions $\varphi_{j,k}^{\left[0\right]}(x) = \chi_{j,k}(x)$, where
$\chi_{j,k}(x) = 1 \Leftrightarrow x \in [x_{j,k},x_{j,k+1})$ and
$\chi_{j,k}(x) = 0$ otherwise. B-splines of order 1 are also known as
B-splines of degree 0.

B-splines of order $\widetilde{p}$, i.e., degree $\widetilde{p}-1$, for
$\widetilde{p}>0$, are defined recursively as
\begin{eqnarray}
\varphi_{j,k}^{\left[\widetilde{p}\right]}(x)
& = &
{x-x_{j,k-\floor{\widetilde{p}/2}} \over
x_{j,k+\ceil{\widetilde{p}/2}-1}-x_{j,k-\floor{\widetilde{p}/2}}}
\varphi_{j,k-1+\mathrm{rem}(\widetilde{p}/2)}^{\left[\widetilde{p}-1\right]}(x)
\nonumber\\
& + &
{x_{j,k+\ceil{\widetilde{p}/2}}-x \over
x_{j,k+\ceil{\widetilde{p}/2}}-x_{j,k-\floor{\widetilde{p}/2}+1}}
\varphi_{j,k+\mathrm{rem}(\widetilde{p}/2)}^{\left[\widetilde{p}-1\right]}(x).
\label{eq:Bsplinerecursion}
\end{eqnarray}
In this equation $\mathrm{rem}(p/q) = p-q\floor{p/q}$ denotes the remainder
from an integer division.
\label{def:Bspline}
\end{definition}
Later on in this paper, the construction through recursion will be replaced by
a construction through refinement. On a finite set of knots, i.e., when $k \in
\{0,\ldots,n_j-1\}$, both constructions are equivalent if we follow the
convention in (\ref{eq:Bsplinerecursion}) that the left and right end points
are multiple knots. More precisely, whenever a knot index in
(\ref{eq:Bsplinerecursion}) is outside $\{0,\ldots,n_j-1\}$, then we take
$x_{j,l} = x_{j,0}$ for $l<0$ and $x_{j,r} = x_{j,n_j-1}$ for $r>n_j-1$.
Definition \ref{eq:Bsplinerecursion}
associates with every knot $x_{j,k}$ a B-spline function
$\varphi_{j,k}^{\left[\widetilde{p}\right]}(x)$. The function turns out to be
centered around the corresponding knot, as can be seen from the following
result.
\begin{thm}(piecewise polynomials on bounded intervals)
For $k \in \{\floor{\widetilde{p}/2},\ldots,n-1-\ceil{\widetilde{p}/2}\}$
the function $\varphi_{j,k}^{\left[\widetilde{p}\right]}(x)$ is zero outside
the interval $[x_{j,l_k},x_{j,r_k})$, where $l_k = k-\floor{\widetilde{p}/2}$
and $r_k = k+\ceil{\widetilde{p}/2}$. Inside this interval,
$\varphi_{j,k}^{\left[\widetilde{p}\right]}(x)$ is a polynomial of
degree $\widetilde{p}-1$ between two knots $x_{j,k}$ and $x_{j,k+1}$,
while in the knots, the function and its first $\widetilde{p}-2$ derivatives
are continuous.
\label{thm:piecewisepolynomial}
\end{thm}

\begin{longversion}
The proof follows by induction, using Definition \ref{def:Bspline}.

Theorem \ref{thm:piecewisepolynomial} should be amended for functions
$\varphi_{j,k}^{\left[\widetilde{p}\right]}$ near the boundaries, i.e.,
for $k$ close to 0 or $n_j-1$. The specification is postponed to the moment
where the functions have been redefined using refinement instead of recursion.
We first concentrate on the interior interval, defined by
\begin{equation}
I_j = \{x \in [x_{j,0},x_{j,n_j-1}) | \forall k=0,\ldots,n_j-1:
\varphi_{j,k}^{\left[\widetilde{p}\right]}(x) \neq 0 \Rightarrow
l_k\geq 0 \mbox{ and } r_k\leq n_j-1\}, 
\label{def:interiorinterval}
\end{equation}
with $l_k$ and $r_k$ as defined in Theorem \ref{thm:piecewisepolynomial}.
It is straightforward to check that
\begin{equation}
I_j = \left[x_{j,\widetilde{p}-1},x_{j,n_j-\widetilde{p}}\right].
\label{eq:Ij}
\end{equation}
\end{longversion}

From Theorem \ref{thm:piecewisepolynomial}, it is obvious that the set of
B-splines $\left\{\varphi_{j,k}^{\left[\widetilde{p}\right]}\right\}$
generates piecewise
polynomials. Conversely, it can be verified that any piecewise polynomial on
the interior interval
\begin{shortversion}
\begin{equation}
I_j = \left[x_{j,\widetilde{p}-1},x_{j,n_j-\widetilde{p}}\right]
\label{eq:Ij}
\end{equation}
\end{shortversion}
\begin{longversion}
$I_j$
\end{longversion}
can be decomposed as a linear combination of B-splines.
\begin{longversion}
\begin{thm}(B-spline basis)
Let $f_j(x)$ be a function which is polynomial on each interval
$[x_{j,k},x_{j,k+1})$ for which $\{x_{j,k},x_{j,k+1}\} \subset I_j$ and which
has $p-2$ continuous derivatives in the knots $x_{j,k} \in I_j$. Then there
exist constants $a_{j,k}$ so that for all $x \in I_j$,
\begin{equation}
f_j(x) = \sum_{k=\floor{\widetilde{p}/2}}^{n_j-1-\ceil{\widetilde{p}/2}}
a_{j,k}\, \varphi_{j,k}^{\left[\widetilde{p}\right]}(x)
\label{eq:Bsplinebasis}
\end{equation}
\label{thm:Bsplinebasis}
\end{thm}
\proof See Appendix \ref{app:proof:thm:Bsplinebasis}.
\end{longversion}
\begin{remark}
Most references in literature would adopt the symbol
$N_{j,k}^{\left[\widetilde{p}\right]}(x)$ for a shifted version of this basis,
namely
\[
\varphi_{j,k}^{\left[\widetilde{p}\right]}(x) =
N_{j,k-\floor{\widetilde{p}/2}}^{\left[\widetilde{p}\right]}(x).
\]
The notation with $N_{j,k}^{\left[\widetilde{p}\right]}(x)$ leads to more
elegant expressions for the recursion of B-splines, but it does not correspond
to the common practice in wavelet literature where a basis function
$\varphi_{j,k}(x)$ is centered around the point $x_{j,k}$.
\end{remark}

The forthcoming discussions will use expressions for the derivatives of
spline functions. 
\begin{shortversion}
In particular the function
\end{shortversion}
\begin{longversion}
\begin{lemma}
The derivative of a B-spline is given by
\begin{equation}
{d\over dx}\varphi_{j,k}^{\left[\widetilde{p}\right]}(x)
= (\widetilde{p}-1)\left[{
\varphi_{j,k-1+\mathrm{rem}(\widetilde{p}/2)}^{\left[\widetilde{p}-1\right]}(x)
\over 
x_{j,k+\ceil{\widetilde{p}/2}-1}-x_{j,k-\floor{\widetilde{p}/2}}}
- 
{\varphi_{j,k+\mathrm{rem}(\widetilde{p}/2)}^{\left[\widetilde{p}-1\right]}(x)
\over 
  x_{j,k+\ceil{\widetilde{p}/2}}-x_{j,k-\floor{\widetilde{p}/2}+1}}
\right].
\label{eq:Bsplinederivative}
\end{equation}
\label{lemma:Bsplinederivative}
\end{lemma}
Lemma \ref{lemma:Bsplinederivative} can be proven by induction on
$\widetilde{p}$.
The computations are facilitated by working on the shifted index in
$N_{j,k}^{\left[\widetilde{p}\right]}(x)$.

As a consequence of Lemma \ref{lemma:Bsplinederivative}, a linear combination
of the B-splines
\end{longversion}
\begin{equation}
f_j(x) = \sum_{k\in \ZZ} s_{j,k}^{\left[\widetilde{p}\right]}
\varphi_{j,k}^{\left[\widetilde{p}\right]}(x),
\label{eq:expansionBsplinebasis}
\end{equation}
has a derivative equal to
\begin{equation}
f_j'(x) = (\widetilde{p}-1) \sum_{k\in \ZZ} 
         {s_{j,k}^{\left[\widetilde{p}\right]} - 
          s_{j,k-1}^{\left[\widetilde{p}\right]} \over 
           x_{j,k+\ceil{\widetilde{p}/2}-1}-x_{j,k-\floor{\widetilde{p}/2}}}
           \varphi_{j,k-\widetilde{r}'}^{\left[\widetilde{p}-1\right]}(x),
\label{eq:derivativeinBsplinebasis}
\end{equation}
where $\widetilde{r}' = 1 - \mathrm{rem}(\widetilde{p}/2)$.

In the search for a refinement relation for B-splines, an important role will
be played by the decomposition of the power functions $x^q$ of degrees
$q=0,\ldots,\widetilde{p}-1$. 
\begin{longversion}
For $q=0$, Definition \ref{def:Bspline} allows us to conclude that the basis
functions are normalised so that
\begin{equation}
\forall x \in I_j:
\sum_{k = 0}^{n_j-1} \varphi_{j,k}^{\left[\widetilde{p}\right]}(x) = 1,
\label{eq:partitionunityBsplines}
\end{equation}
for any order $\widetilde{p}$. This property is referred to as the
\emph{partition of unity}.
For general $q<\widetilde{p}$, the expansion of $x^q$ in a B-spline basis can
be established according to the following result, which is closely related to
Marsden's identity \citep{lee96:marsden}.
\end{longversion}
\begin{shortversion}
Their expansion in a spline basis can be established according to the 
following result, which is closely related to Marsden's identity
\citep{lee96:marsden}.
\end{shortversion}
\begin{thm} (power coefficients)
For $q=0,1,\ldots,\widetilde{p}-1$, there exist coefficients
$\widetilde{x}_{j,k}^{\left[\widetilde{p},q\right]}$, so that for $x \in I_j$,
\begin{equation}
x^q = \sum_{k\in \ZZ} \widetilde{x}_{j,k}^{\left[\widetilde{p},q\right]}
\varphi_{j,k}^{\left[\widetilde{p}\right]}(x).
\label{eq:xpowerqBsplinedecomposition}
\end{equation}
\begin{enumerate}
\item
\begin{shortversion}
For $q=0$, these coefficients are one. 
In the wavelet literature, this property is refered to as the partition of
unity.
\end{shortversion}
\begin{longversion}
For $q=0$, these coefficients are one, according to the partition of unity.
\end{longversion}
\item
For $q=1,\ldots,\widetilde{p}-1$, the coefficients can be found by the 
recursion
\begin{equation}
\widetilde{x}_{j,k}^{\left[\widetilde{p},q\right]} =
\widetilde{x}_{j,k-1}^{\left[\widetilde{p},q\right]} +
{q \over \widetilde{p}-1}\,
\widetilde{x}_{j,k-\widetilde{r}'}^{\left[\widetilde{p}-1,q-1\right]}\,
\left(x_{j,k+\ceil{\widetilde{p}/2}-1}-x_{j,k-\floor{\widetilde{p}/2}}\right),
\label{eq:polynomialcoefsBspline}
\end{equation}
where $\widetilde{r}' = 1 - \mathrm{rem}(\widetilde{p}/2)$, as in
(\ref{eq:derivativeinBsplinebasis}).
\item
In particular, for $q=1$, these coefficients satisfy
\begin{equation}
\widetilde{x}_{j,k}^{\left[\widetilde{p},1\right]}
=
{1 \over \widetilde{p}-1}
\sum_{i=1-\floor{\widetilde{p}/2}}^{\ceil{\widetilde{p}/2}-1} x_{j,k+i}.
\label{eq:linearcoefsBspline}
\end{equation}
\item
For $q = \widetilde{p}-1$, this becomes
\begin{equation}
\widetilde{x}_{j,k}^{\left[\widetilde{p},\widetilde{p}-1\right]}
=
\prod_{i=1-\floor{\widetilde{p}/2}}^{\ceil{\widetilde{p}/2}-1} x_{j,k+i}.
\label{eq:x^(p-1)coefsBspline}
\end{equation}
\end{enumerate}
\label{thm:x^qcoefBspline}
\end{thm}
\begin{longversion}
\proof See Appendix \ref{app:proof:thm:x^qcoefBspline}.
\end{longversion}

\begin{remark}
Section \ref{subsec:designliftingsteps} will be based on a general reading
of (\ref{eq:polynomialcoefsBspline}), which describes the transition
from $\widetilde{x}_{j,k-1}^{\left[\widetilde{p},q\right]}$ to
$\widetilde{x}_{j,k}^{\left[\widetilde{p},q\right]}$. Since the results in
Theorem \ref{thm:x^qcoefBspline} do not depend on the ordering of the knots,
the same formula as in (\ref{eq:polynomialcoefsBspline}) can also be used to
recompute the coefficients
$\widetilde{x}_{j,k-1}^{\left[\widetilde{p},q\right]}$ when one
knot is taken out from the grid and replaced by another. 
In (\ref{eq:polynomialcoefsBspline}) the knot $x_{j,k-\floor{\widetilde{p}/2}}$
is replaced by $x_{j,k+\ceil{\widetilde{p}/2}-1}$, while the factor
$\widetilde{x}_{j,k-\widetilde{r}'}^{\left[\widetilde{p}-1,q-1\right]}$ depends
only on knots that are left untouched.

Expression (\ref{eq:polynomialcoefsBspline}) is an example of a formula that
is simplified by working in the shifted basis 
$N_{j,k}^{\left[\widetilde{p}\right]}(x)$.
Putting $t_{j,k} = x_{j,k+\floor{\widetilde{p}/2}}$, and
\(
\widetilde{t}_{j,k}^{\left[\widetilde{p},q\right]} =
\widetilde{x}_{j,k+\floor{\widetilde{p}/2}}^{\left[\widetilde{p},q\right]},
\)
we get
\begin{equation}
\widetilde{t}_{j,k}^{\left[\widetilde{p},q\right]} =
\widetilde{t}_{j,k-1}^{\left[\widetilde{p},q\right]} +
{q \over \widetilde{p}-1}\,
\widetilde{t}_{j,k}^{\left[\widetilde{p}-1,q-1\right]}\,
\left(t_{j,k+\widetilde{p}-1}-t_{j,k}\right).
\label{eq:polynomialcoefsBsplineN}
\end{equation}
\end{remark}

The following theorem states that no other basis reproduces polynomials with
functions that have a support between $\widetilde{p}+1$ knots.
\begin{thm}(uniqueness by power coefficients)
Let $x_{j,k}$ for $k = 0,\ldots,n_j-1$ be the knots $x_{j,k}$ at level $j$, 
and let $\varphi_{j,k}^{\left[\widetilde{p}\right]}(x)$ be a set of basis
functions associated to these knots. If the support of
$\varphi_{j,k}^{\left[\widetilde{p}\right]}(x)$ equals
$S_{j,k} = [x_{j,l_k},x_{j,r_k})$,
with $l_k$ and $r_k$ as in Theorem \ref{thm:piecewisepolynomial},
and if the coefficients $\widetilde{x}_{j,k}^{\left[\widetilde{p},q\right]}$
in the decompositions (\ref{eq:xpowerqBsplinedecomposition}),
for all $x \in I_j$ and for all $q \in \{0,1,\ldots,\widetilde{p}-1\}$,
are given by the values in Theorem \ref{thm:x^qcoefBspline},
then $\varphi_{j,k}^{\left[\widetilde{p}\right]}(x)$ must be the B-splines of
order $\widetilde{p}$ defined on the given knots.
\label{thm:Bsplinesfrompolynomialreproduction}
\end{thm}
\begin{longversion}
\proof See Appendix \ref{app:proof:thm:Bsplinesfrompolynomialreproduction}.
\end{longversion}
Theorem \ref{thm:Bsplinesfrompolynomialreproduction} motivates the use of the
power function coefficients
$\widetilde{x}_{j,m}^{\left[\widetilde{p},q\right]}$ as a defining property
in the design of a refinement scheme for B-splines.

\begin{shortversion}
\subsection{B-spline refinement schemes and factorisation}
\label{subsec:Bspline2scales}
\end{shortversion}

\begin{longversion}
\subsection{B-spline refinement schemes}
\label{subsec:Bspline2scales}
\end{longversion}

In this section, we consider B-spline functions at different scales, i.e., with
different indices $j$. 
\begin{shortversion}
It can be verified that a construction through refinement must exist
\citep[Eq.(26), page 497]{qu92:subdivisionBsplines,daubechies01:commutation}.
i.e., the B-splines on any set of knots must satisfy a refinement or two-scale
equation, whose general form is
\end{shortversion}
\begin{longversion}
The first result states that a construction through refinement must exist, see
\citet[(26), page 497]{qu92:subdivisionBsplines,daubechies01:commutation}.
\begin{thm} (existence of B-spline refinement)
On a nested multilevel grid, B-spline basis functions at scale $j$ are
refinable, i.e., there exists a refinement matrix
$H_j^{\left[\widetilde{p}\right]}$ so that, for all $x \in I_j$,
\begin{equation}
\varphi_{j,k}^{\left[\widetilde{p}\right]}(x)
=
\sum_{\ell=0}^{n_{j+1}} H_{j,\ell,k}^{\left[\widetilde{p}\right]}\,
\varphi_{j+1,\ell}^{\left[\widetilde{p}\right]}(x).
\label{eq:2scaleBsplines}
\end{equation}
\end{thm}

\proof
This is an immediate consequence of Theorems \ref{thm:piecewisepolynomial}
and \ref{thm:Bsplinebasis}.
A B-spline on a grid at level $j$ is a piecewise polynomial with knots in
$x_{j,k}$. Since these knots are also knots at scale $j+1$, the B-spline is
also a piecewise polynomial at scale $j+1$, and so it can be written as a
linear combination of the B-spline basis at that scale.
\qed\\
Equation (\ref{eq:2scaleBsplines}) is an instance of a two-scale equation,
also known as a refinement equation. The general form of the refinement
equation, without superscripts for the order of the B-splines, is
\end{longversion}
\begin{equation}
\varphi_{j,k}(x)
=
\sum_{\ell=0}^{n_{j+1}} H_{j,\ell,k} \varphi_{j+1,\ell}(x).
\label{eq:2scale}
\end{equation}
The refinement equation can also be condensed into a matrix form
\begin{equation}
\Phi_j(x) = \Phi_{j+1}(x)H_j,
\label{eq:2scalematrix}
\end{equation}
where
\begin{equation}
\Phi_j(x) = [\varphi_{j,0}(x) \, \varphi_{j,1}(x) \, \ldots
\,\varphi_{j,n_j-1}(x)]
\end{equation}
is a row of scaling basis functions.

\begin{shortversion}
A refinement matrix $H_j$ is typically band-limited, and this holds for
B-spline refinement as well.
\begin{lemma}
Let $H_j^{\left[\widetilde{p}\right]}$ be the refinement matrix for the
B-splines as scaling functions in (\ref{eq:2scalematrix}), then the columns of
$H_j^{\left[\widetilde{p}\right]}$ can have at most $\widetilde{p}+1$ nonzero
elements. In particular 
\begin{equation}
H_{j,\ell,k}^{\left[\widetilde{p}\right]} \neq 0
\Rightarrow 
2k-\floor{\widetilde{p}/2} \leq \ell \leq 2k+\ceil{\widetilde{p}/2}.
\end{equation}
\end{lemma}
\end{shortversion}
\begin{longversion}
In the first instance, a refinement equation should be read as the definition
of the scaling functions from the refinement matrices $H_j$. A numerical
solution of the equation thus allows us to evaluate the scaling functions, in
particular the B-spline functions.
Secondly, the refinement equation will be the basis for the
construction of B-spline wavelets, as explained in Section
\ref{subsec:Bsplinewavelets}.
This motivates the search for the spline refinement matrix
$H_j^{\left[\widetilde{p}\right]}$.

A refinement matrix $H_j$ is often band-limited, as follows from the next
lemma, valid for general refinable scaling functions.
\begin{lemma}(band-limited refinement matrices)
Let $\varphi_{j,k}(x)$, for $j=0,\ldots,n_j-1$ be a set of scaling functions,
with refinement equation (\ref{eq:2scale}).
Let $S_{j,k}$ denote the support of $\varphi_{j,k}(x)$, then
an entry $H_{j,\ell,k}$ may be different
from zero only if $S_{j+1,\ell} \subset S_{j,k}$.
\label{lemma:Hbandlimited}
\end{lemma}
\proof
Suppose $H_{j,\ell,k}$ would be nonzero for a fine scaling function
outside the support of the coarse scaling function, then obviously, that
coarse scaling function would take a nonzero value outside its support.
\qed

For the B-spline basis, Theorem \ref{lemma:Hbandlimited} translates as follows.
\begin{corollary}
The columns of the matrix $H_j^{\left[\widetilde{p}\right]}$ in
(\ref{eq:2scaleBsplines}) can have at most $\widetilde{p}+1$ nonzero
elements. In particular
\begin{equation}
H_{j,\ell,k}^{\left[\widetilde{p}\right]} \neq 0
\Rightarrow 
2k-\floor{\widetilde{p}/2} \leq \ell \leq 2k+\ceil{\widetilde{p}/2}.
\end{equation}
\end{corollary}
\proof
The support of a B-spline $\varphi_{j,k}^{\left[\widetilde{p}\right]}(x)$ is
$S_{j,k} = [x_{j,k-\floor{\widetilde{p}/2}},x_{j,k+\ceil{\widetilde{p}/2}}]$.
We have
\begin{eqnarray*}
S_{j+1,\ell} \subset S_{j,k} 
& \Leftrightarrow &
[x_{j+1,\ell-\floor{\widetilde{p}/2}},x_{j+1,\ell+\ceil{\widetilde{p}/2}}]
\subset
[x_{j,k-\floor{\widetilde{p}/2}},x_{j,k+\ceil{\widetilde{p}/2}}]
\\
& \Leftrightarrow &
[x_{j+1,\ell-\floor{\widetilde{p}/2}},x_{j+1,\ell+\ceil{\widetilde{p}/2}}]
\subset
[x_{j+1,2k-2\floor{\widetilde{p}/2}},x_{j+1,2k+2\ceil{\widetilde{p}/2}}]
\\
& \Leftrightarrow &
\left\{\begin{array}{l} \ell-\floor{\widetilde{p}/2}\geq
2k-2\floor{\widetilde{p}/2}
\mbox{ and}\\
                        \ell+\ceil{\widetilde{p}/2} \leq
2k+2\ceil{\widetilde{p}/2}
       \end{array}\right.
\\
& \Leftrightarrow &
\ell \in \{2k-\floor{\widetilde{p}/2},\ldots,2k+\ceil{\widetilde{p}/2}\}.
\end{eqnarray*}
As $\#\{2k-\floor{\widetilde{p}/2},\ldots,2k+\ceil{\widetilde{p}/2}\} =
\ceil{\widetilde{p}/2} + \floor{\widetilde{p}/2} +1 =
\widetilde{p}+1$, the number of nonzero elements in column $k$ is bounded by
$\widetilde{p}+1$.
\qed\\
The maximum number of nonzero elements in each column is known as the bandwidth
of the matrix, see Definition \ref{def:bandwidth} below.

\subsection{Factorisation of the refinement matrix}
\label{subsec:factoring2scaleeq}
\end{longversion}

The objective is of course to identify the nonzero entries of
$H_j^{\left[\widetilde{p}\right]}$. The strategy followed in this paper is
based on a factorisation that can be applied to any band-limited refinement
matrix $H_j$.
The factorisation starts from a partition of the rows of
the matrix into an even and an odd subset, leading to submatrices
$H_{j,e}$ and $H_{j,o}$
and so that the refinement equation (\ref{eq:2scalematrix}) can be written as
\begin{equation}
\Phi_j(x) = \Phi_{j+1,e}(x)H_{j,e}+ \Phi_{j+1,o}(x)H_{j,o}.
\label{eq:2scalePhiHeHo}
\end{equation}
\begin{remark}
The matrix $H_{j,e}$ is a squared matrix, because in the nested refinement, the
even subset of knots at scale $j+1$ are exactly the knots at scale $j$.
\end{remark}

These submatrices can then be factored in an iterative way, alternating
between two sorts of factorisation steps. The alternation is the matrix
equivalent of Euclid's algorithm for finding the greatest common
divider \citep{daubechies98:factorlifting}. The superscript $[s]$ in the
theorem refers to the factorisation step. It
should not be confused with the superscript $\left[\widetilde{p}\right]$
referring to the order of a B-spline. In Section
\ref{subsec:designliftingsteps}, both $s$ and
$\widetilde{p}$ will appear in single superscripts.
The subsequent result uses the following definition for bandwidth of a
rectangular matrix.
\begin{definition} (bandwidth of a refinement matrix)
Let $A$ be an $m \times n$ matrix, where in each column $j = 1,\ldots,n$, there
exists a row $i_1(j)$ so that
\begin{equation}
A_{i,j} \neq 0 \Rightarrow 
i_1(j) \leq i \leq i_1(j)+b-1,
\label{eq:bandwidth}
\end{equation}
with $b$ independent from $j$,
then the bandwidth of this matrix is $b$.
\label{def:bandwidth}
\end{definition}
\begin{thm}(factorisation into lifting steps)
Given a refinement matrix $H_j^{[s]}$ and the submatrices 
$H_{j,e}^{[s]}$ and $H_{j,o}^{[s]}$
containing its even and odd rows respectively. If $H_{j,e}^{[s]}$ has a larger
bandwidth than $H_{j,o}^{[s]}$, then we can always find a lower
\emph{bidiagonal} matrix $U_j^{[s+1]}$ and a matrix $H_{j,e}^{[s+1]}$ with
smaller bandwidth than that of $H_{j,o}^{[s]}$ so that
\begin{equation}
H_{j,e}^{[s]} = H_{j,e}^{[s+1]} - U_j^{[s+1]} H_{j,o}^{[s]},
\label{eq:primallifting}
\end{equation}
If $H_{j,o}^{[s]}$ has a larger bandwidth than $H_{j,e}^{[s]}$, then we can
always find an upper \emph{bidiagonal} matrix $P_j^{[s+1]}$ and a matrix
$H_{j,o}^{[s+1]}$ with smaller bandwidth than
that of $H_{j,e}^{[s]}$ so that
\begin{equation}
H_{j,o}^{[s]} = H_{j,o}^{[s+1]} + P_j^{[s+1]} H_{j,e}^{[s]},
\label{eq:duallifting}
\end{equation}
\label{thm:factorisation}
If $H_{j,e}^{[s]}$ and $H_{j,o}^{[s]}$ have the same bandwidth, then both
(\ref{eq:primallifting}) and (\ref{eq:duallifting}) are possible.
\end{thm}
\proof See Appendix \ref{app:proof:thm:factorisation}.

The matrices $P_j^{[s+1]}$ are known as dual lifting steps or prediction steps,
where both terms refer to an interpretation beyond the scope of this paper
\citep{sweldens97:lifting2nd}.
As a matter of fact, the interpretation as a prediction is even not applicable
in the context of this paper. The matrices $U_j^{[s+1]}$ are primal lifting
steps or update steps.

The next sections explain how the factorisation into lifting steps can be used
in the design of a multiscale decomposition of a B-spline basis on irregular
knots.


\section{Non-equispaced B-spline wavelet transforms}


\subsection{Main contribution: the design of B-spline lifting steps}
\label{subsec:designliftingsteps}

\begin{revision}
Theorem \ref{thm:Bsplinesfrompolynomialreproduction} allows us to develop a
lifting scheme for B-splines resting on an analysis of power function
coefficients only.
\end{revision}
We will impose that if $s_{j+1,k} =
\widetilde{x}_{j+1,k}^{\left[\widetilde{p},q\right]}$, then
$d_{j,k}$ must be zero, while
$s_{j,k} = \widetilde{x}_{j,k}^{\left[\widetilde{p},q\right]}$.
In principle, the development of this condition can proceed directly on the
band matrix $H_j^{\left[\widetilde{p}\right]}$, without having to use the
lifting factorisation. Solving the corresponding linear system seems, however,
not to yield a simple formula.

The lifting factorisation of $H_j^{\left[\widetilde{p}\right]}$ is found in a
relatively straightforward way, because all lifting steps are essentially the
same sort of operation, described in the following proposition.
The successive lifting steps differ from each other only in the subsets of
knots that are involved in the operation. As a consequence, we formulate the
operation for an arbitrary subset of knots. The elements of the subsets are
denoted by $t_{j,i}$. Obviously, all $t_{j,i}$ coincide with a knot from
the $x_k$, but the $t_{j,i}$ may be unordered and the subset may contain
coinciding values $t_{j,i'} = t_{j,i}$. The proposition constructs a recursion
(\ref{eq:polynomialcoefsBsplineN}) on the selected knots $t_{j,i}$. The
recursion (\ref{eq:polynomialcoefsBsplineN}) has originally been stated for
$t_{j,i}$ that are shifted sorted knots, but nothing prevents us from using it
for more general sets of $t_{j,i}$.
\begin{prop} (lifting step for B-spline power coefficients)
Given an order $\widetilde{p}$,
let $K_{j,k} = \{t_{j,k},\ldots,t_{j,k+\widetilde{p}}\}$ an unsorted set of
knots.
Define the coefficient $\widetilde{t}_{j,k}^{\left[\widetilde{p},q\right]}$
by the recursion in (\ref{eq:polynomialcoefsBsplineN}), using the knots
$\{t_{j,k+1},\ldots,t_{j,k+\widetilde{p}-1}\}$. Define the coefficients
$\widetilde{t}_{j,k-1}^{\left[\widetilde{p},q\right]}$ and
$\widetilde{t}_{j,k+1}^{\left[\widetilde{p},q\right]}$ by the same recursion,
but using $\{t_{j,k},\ldots,t_{j,k+\widetilde{p}-2}\}$ and
$\{t_{j,k+2},\ldots,t_{j,k+\widetilde{p}}\}$ respectively.

Define the lifting parameters
\begin{equation}
L_{j,k,k-1} = {t_{j,k+\widetilde{p}}-t_{j,k+\widetilde{p}-1} \over
                   t_{j,k+\widetilde{p}}-t_{j,k}}
\,\,\,\mbox{ and }\,\,\,
L_{j,k,k+1} = {t_{j,k+1}-t_{j,k} \over t_{j,k+\widetilde{p}}-t_{j,k}}.
\label{eq:liftingpowercoef}
\end{equation}
Then the lifted coefficient
\(
\widetilde{t}_{j,k}^{\left[\widetilde{p},q,L1\right]} =
\widetilde{t}_{j,k}^{\left[\widetilde{p},q\right]}
-
L_{j,k,k-1}
\widetilde{t}_{j,k-1}^{\left[\widetilde{p},q\right]}
-
L_{j,k,k+1}
\widetilde{t}_{j,k+1}^{\left[\widetilde{p},q\right]}
\label{eq:liftedpowercoef}
\)
equals, for all values of $q=0,\ldots,\widetilde{p}-1$,
\begin{equation}
\widetilde{t}_{j,k}^{\left[\widetilde{p},q,L1\right]} =
{t_{j,k+\widetilde{p}-1}-t_{j,k+1} \over t_{j,k+\widetilde{p}}-t_{j,k}}
\,
\widetilde{t}_{j,k}^{\left[\widetilde{p},q,L1b\right]},
\label{eq:liftingoutandin}
\end{equation}
where $\widetilde{t}_{j,k}^{\left[\widetilde{p},q,L1b\right]}$ is the
coefficient defined by the recursion in (\ref{eq:polynomialcoefsBsplineN}),
using the knots\\
$\{t_{j,k},t_{j,k+2},\ldots,t_{j,k+\widetilde{p}-2},t_{j,k+\widetilde{p}}\}$.
\label{prop:liftedcoefficient}
\end{prop}
In other words, a lifting step with only two parameters, has a common effect on
all $\widetilde{p}$ power coefficients: it takes out the knots $t_{j,k+1}$ and
$t_{j,k+\widetilde{p}-1}$ to replace them by $t_{j,k}$ and
$t_{j,k+\widetilde{p}}$.
The lifting step has thus a similar effect as the recursive formula
(\ref{eq:polynomialcoefsBsplineN}). The difference between lifting and
recursion is that the recursion uses coefficients of different power functions
in B-splines of different degree, while lifting is based on coefficients of a
single function in a single basis. Moreover, the lifting formula is the same
for all power functions in that basis.

The proof of this theorem is based on the recursion of
(\ref{eq:polynomialcoefsBsplineN}). As it is rather technical, it can be found
in Appendix \ref{app:proof:prop:liftedcoefficient}.

The lifting operation presented in Theorem \ref{prop:liftedcoefficient} lies at
the heart of the linear transform that maps fine scale power coefficients
$\widetilde{x}_{j+1,k}^{\left[\widetilde{p},q\right]}$ onto coarse scale
versions $\widetilde{x}_{j,k}^{\left[\widetilde{p},q\right]}$ along with zero
detail coefficients . We adopt the symbol
$\widetilde{x}_{j,k}^{\left[\widetilde{p},q\right]}$ in contrast to
$\widetilde{t}_{j,k}^{\left[\widetilde{p},q\right]}$ to emphasise that we are
working now on the sorted knots $x_{j,k}$.
The update lifting steps $U_j^{[s]}$ take care of coarse scaling
coefficients, starting from the even fine scale coefficients. Each
update step takes out two odd indexed, fine scale knots $x_{j+1,2k\pm (2m+1)}$
from the intermediate scaling coefficients and adds two coarse scale
knots $x_{j,k\pm t} = x_{j+1,2k\pm 2t}$ that are outside the fine scale
range. The indices $m$ and $t$ depend on the lifting step $s$, as developed in
Proposition \ref{prop:liftingschemeBsplines}.
In a similar way, the prediction lifting steps $P_j^{[s]}$ take care of the
detail coefficients, operating on the odd fine scaling coefficients. A
prediction step takes out two odd indexed, fine scale knots $x_{j+1,2k\pm
2m+1}$ from the definition of the intermediate coefficient, replacing it by
two new coarse scale knots.

The final prediction step
has a special role. It is supposed to take out the
last remaining odd knot twice. That is, $m = 0$, so that $x_{j+1,2k+2m+1} =
x_{j+1,2k-2m+1}$. In terms of the unsorted knots $t_{j,k}$ in
Proposition \ref{prop:liftedcoefficient}, the
knots are numbered so that $t_{j,k+\widetilde{p}-1}=t_{j,k+1}$. The outcome
of the lifted power coefficient in (\ref{eq:liftedpowercoef}) is then zero, as
requested. The preceeding prediction steps should therefore be such that one
odd knot is left over for the final step. Depending on the number of odd
knots in the beginning, this may imply that the first prediction step takes
out only one odd knot. 
\begin{revision}
In that case, the first prediction step is a diagonal instead of a bidiagonal
matrix. Whether a prediction matrix has one or two nonzero (off-)diagonals is
controlled by the variable $\widetilde{t}$, defined in Proposition
\ref{prop:liftingschemeBsplines}.
Similar considerations hold for the update steps, for which the variable
$\widetilde{u}$ controls the number of nonzero diagonals.
\end{revision}
All together, we arrive at the following lifting scheme for B-splines.
\begin{prop} (Main result: lifting scheme for B-splines)
Suppose that $\vec{s}_{j+1,k}$ are scaling coefficients at scale $j+1$ defined
on the knots $\vec{x}_{j+1,k}$. Then consider a lifting scheme with $u$ update
steps and $u+r$ prediction steps, where integer $u$ and boolean value $r$ are
given by
\(
u = \floor{(\widetilde{p}+1)/4},
\)
and
\(
r = \ceil{\widetilde{p}/2} - 2u,
\)
for a given integer $\widetilde{p}$. Furthermore, let $\widetilde{r} =
\widetilde{p}-2\floor{\widetilde{p}/2}$ be a boolean indicating the parity of
$\widetilde{p}$.
For every $s \in \{1-r,\ldots,u\}$, define the values $m = u-s$, $t=r+2s-1$,
$\widetilde{u} = \widetilde{r}\cdot\min(r+s,2)$, 
$\widetilde{t} = \widetilde{r}\cdot\min(1+s,2)$. 

Then construct the following lifting scheme.
First, define the even and odd factors, for $s \in \{1,\ldots,u\}$, and for $s
\in \{1-r,\ldots,u\}$, respectively,
\begin{eqnarray}
c_{j+1,k}^{[0]} & = & 1
\label{eq:factorsliftingBsplines0}
\\
c_{j+1,2k}^{[s]}
& = &
c_{j+1,2k}^{[s-1]} \cdot
{x_{j+1,2k+2m+1}-x_{j+1,2k-2m-1} \over
 x_{j+1,2k+2m+2t}-x_{j+1,2k-2m-2t+\widetilde{u}}}
\label{eq:factorsliftingBsplinesE}
\\
c_{j+1,2k+1}^{[r+s]}
& = &
c_{j+1,2k+1}^{[r+s-1]} \cdot
{x_{j+1,2k+2m+1}-x_{j+1,2k-2m+1} \over
 x_{j+1,2k+2m+2t+2}-x_{j+1,2k-2m-2t+\widetilde{t}}}.
\label{eq:factorsliftingBsplinesO}
\end{eqnarray}
Second, for $s \in \{1,\ldots,u\}$, define an update matrix $U_j^{[s]}$ as a
lower bidiagonal matrix with entries
\begin{eqnarray}
U_{j,k,k}^{[s]} & = & -{c_{j+1,2k}^{[s-1]} \over c_{j+1,2k+1}^{[r+s-1]}}
\cdot
{x_{j+1,2k-2m-1} - x_{j+1,2k-2m-2t+\widetilde{u}} \over
 x_{j+1,2k+2m+2t} - x_{j+1,2k-2m-2t+\widetilde{u}}}
\label{eq:Ujkk}
\\
U_{j,k,k-1}^{[s]} & = & -{c_{j+1,2k}^{[s-1]} \over c_{j+1,2k-1}^{[r+s-1]}}
\cdot
{x_{j+1,2k+2m+2t} - x_{j+1,2k+2m+1} \over
 x_{j+1,2k+2m+2t} - x_{j+1,2k-2m-2t+\widetilde{u}}}.
\label{eq:Ujkk1}
\end{eqnarray}
\begin{revision}
Set the corresponding lifted even coefficients
\begin{equation}
s_{j+1,2k}^{[s]} = s_{j+1,2k}^{[s-1]}
                    + U_{j,k,k-1}^{[s]} s_{j+1,2k-1}^{[r+s-1]}
                    + U_{j,k,k}^{[s]} s_{j+1,2k+1}^{[r+s-1]}.
\end{equation}
Third, for $s \in \{1-r,\ldots,u\}$, define a prediction matrix $P_j^{[r+s]}$
as an upper bidiagonal matrix with entries
\end{revision}
\begin{eqnarray}
P_{j,k,k}^{[r+s]} & = & {c_{j+1,2k+1}^{[r+s-1]} \over c_{j+1,2k}^{[s]}}
\cdot
{x_{j+1,2k+2m+2t+2} - x_{j+1,2k+2m+1} \over
 x_{j+1,2k+2m+2t+2} - x_{j+1,2k-2m-2t+\widetilde{t}}}
\label{eq:Pjkk}
\\
P_{j,k,k+1}^{[r+s]} & = & {c_{j+1,2k+1}^{[r+s-1]} \over c_{j+1,2k+2}^{[s]}}
\cdot
{x_{j+1,2k-2m+1} - x_{j+1,2k-2m-2t+\widetilde{t}} \over
 x_{j+1,2k+2m+2t+2} - x_{j+1,2k-2m-2t+\widetilde{t}}}.
\label{eq:Pjkk1}
\end{eqnarray}
Set the corresponding lifted odd coefficients
\begin{equation}
s_{j+1,2k+1}^{[r+s]} = s_{j+1,2k+1}^{[r+s-1]}
                       - P_{j,k,k}^{[r+s]} s_{j+1,2k}^{[s]}
                       - P_{j,k,k+1}^{[r+s]} s_{j+1,2k+2}^{[s]}.
\end{equation}
Finally, define the diagonal rescaling matrix $D_j$ as
\begin{equation}
D_{j,k,k} = c_{j+1,2k}^{[u]},
\end{equation}
and the scaling and detail coefficients at scale $j$ as
\begin{eqnarray}
s_{j,k} & = & D_{j,k,k}^{-1} s_{j+1,2k}^{[u]}\\
d_{j,k} & = & s_{j+1,2k+1}^{[r+u]}.
\end{eqnarray}
Then, the power coefficients in a B-spline basis at scale $j+1$, defined in
Theorem \ref{thm:x^qcoefBspline}, and denoted as
$s_{j+1,k} = \widetilde{x}_{j+1,k}^{\left[\widetilde{p},q\right]}$,
are transformed by this lifting schemes into coarse scaling coefficients
$s_{j,k} = \widetilde{x}_{j,k}^{\left[\widetilde{p},q\right]}$ plus detail
coefficients $d_{j,k} = 0$. Consequently, by Theorem
\ref{thm:Bsplinesfrompolynomialreproduction}, the refinement equation
associated to this lifting scheme has the B-splines of order $\widetilde{p}$ on
the non-equidistant knots $x_{j,k}$ as its solution.
\label{prop:liftingschemeBsplines}
\end{prop}

\begin{remark}
The lifting scheme for B-splines should thus be understood as a gradual
coarsening of fine scale representation of polynomials. This is unlike some
other lifting constructions, where lifting steps are designed as a way to
improve, i.e., ``lift higher'' existing wavelet transforms, by gradually adding
more properties. In the B-spline case, one could for instance, try to lift a
linear B-spline into a cubic B-spline. Such a construction is, however,
impossible with a single lifting step.
\end{remark}

\begin{revision}
\subsection{Examples of B-spline lifting schemes}
\label{subsec:examplesBsplineliftingschemes}

This section develops concrete examples of Proposition
\ref{prop:liftingschemeBsplines}.

\subsubsection{The Haar scaling functions}

The simplest case of a B-spline basis is that of B-splines of order one, i.e.,
degree zero. The basis functions are characteristic functions on the intervals
between two knots. This is the Haar scaling basis defined on these knots.

As $\widetilde{p}=1$, we find $u=0$, meaning that the lifting scheme has zero
update steps. There will be one prediction step because $r = 1$.
This single prediction step is defined by (\ref{eq:Pjkk}) and (\ref{eq:Pjkk1})
with indices $s=0$, $t=0$, $m=0$, and $r+s-1=0$. 
From (\ref{eq:factorsliftingBsplines0}), it follows that all factors 
$c_{j+1,\ell}^{[0]}$ are equal to one in this prediction. 
We also find $\widetilde{r} = 1$, and so, $\widetilde{t} = 1$,
meaning that the prediction will be a diagonal matrix. This is confirmed by
substitution of all the indices in (\ref{eq:Pjkk}) and (\ref{eq:Pjkk1}).
We find that $P_{j,k,k}^{[1]} = 1$ and $P_{j,k,k+1}^{[1]} = 0$.
We also find that $D_{j,k,k} = c_{j+1,2k}^{[0]} = 1$.

This lifting scheme defines the refinement equation and hence the Haar scaling
functions, in a way developed in Section \ref{subsec:multiscaletransform}.
It does not yet fix the Haar wavelet basis. There are several options for these
basis functions, including the classical Haar basis $\psi_{j,k}(x) =
\varphi_{j+1,2k+1}(x)-\varphi_{j+1,2k}(x)$, but also the Unbalanced Haar basis
\citep{girardi97:unbalanced}. Each option can be realised by one additional
update step, as explained in Section \ref{subsec:Bsplinewavelets}.

\subsubsection{The linear B-spline scaling functions}

Linear splines reproduce constant and linear functions, hence
$\widetilde{p}=2$. We find that $u=0$, so the lifting scheme for the refinement
equation has again no update step. As before, there is one prediction step,
and also the indices $s=0$, $t=0$, $m=0$, and $r+s-1=0$ remain the same as in
the Haar case, again leading to the conclusion that the factors
$c_{j+1,\ell}^{[0]}$ are equal to one. In contrast to the Haar case, we now
have $\widetilde{r} = 0$, and from there, $\widetilde{t} = 0$. The effect of
this is that the prediction is now a bidiagonal matrix, with elements
$P_{j,k,k}^{[1]} = (x_{j+1,2k+2}-x_{j+1,2k+1})/(x_{j+1,2k+2}-x_{j+1,2k})$ and
$P_{j,k,k+1}^{[1]} = (x_{j+1,2k+1}-x_{j+1,2k})/(x_{j+1,2k+2}-x_{j+1,2k})$.
We find again $D_{j,k,k} = c_{j+1,2k}^{[0]} = 1$.

This matrix can be interpreted as a linear interpolation in the odd covariate
values. Therefore, the constant and linear splines have refinement schemes
consisting of, respectively, constant and linear polynomial interpolation
in a single prediction step.
Higher order splines cannot be associated with higher order polynomial
interpolation, as illustrated by the next example.

\subsubsection{The cubic B-spline scaling functions}

For cubic B-splines, we set $\widetilde{p}=4$, from which it follows that
$u=1$, meaning that we have one update step. There will be one prediction step,
as $r = 0$, so the lifting scheme starts with the update.
The index $s \in \{1-r,\ldots,u\} = \{1\}$ takes only one value, and so do the
indices $t=1$, $m=0$, and $r+s-1=0$. Furthermore $\widetilde{r} = 0$, from
which it follows that $\widetilde{t} = 0 = \widetilde{u}$, meaning that both
the update and the prediction are bidiagonal matrices. For the update, we find
$U_{j,k,k}^{[1]} = -(x_{j+1,2k-1}-x_{j+1,2k-2})/(x_{j+1,2k+2}-x_{j+1,2k-2})$ 
and
$U_{j,k,k-1}^{[1]} = -(x_{j+1,2k+2}-x_{j+1,2k+1})/(x_{j+1,2k+2}-x_{j+1,2k-2})$.
The subsequent prediction is given by
$P_{j,k,k}^{[1]} = (x_{j+1,2k+4}-x_{j+1,2k+1})/(x_{j+1,2k+4}-x_{j+1,2k-2})$ and
$P_{j,k,k+1}^{[1]} = (x_{j+1,2k+1}-x_{j+1,2k-2})/(x_{j+1,2k+4}-x_{j+1,2k-2})$.
The diagonal elements of the final rescaling become
\(
D_{j,k,k} = c_{j+1,2k}^{[1]} =
(x_{j+1,2k+1}-x_{j+1,2k-1})/(x_{j+1,2k+2}-x_{j+1,2k-2}).
\)

\end{revision}

\subsection{From the factorisation to a multiscale transform}
\label{subsec:multiscaletransform}

\begin{revision}
The two-scale equation (\ref{eq:2scalematrix}) containing the matrix $H_j$
in Section \ref{subsec:Bspline2scales} concentrated on the refinement of
coarse scale functions $f_j(x)$ or coarse scale basis functions $\Phi_j(x)$. 
The argument used for the design of $H_j$ in Proposition
\ref{prop:liftingschemeBsplines} went the other way, starting from the finer
scale $j+1$, and imposing that fine scale power coefficients are projected onto
coarse scale power projections. This section will reassemble the matrix $H_j$
from the lifting factorisation that resulted from Proposition
\ref{prop:liftingschemeBsplines}. In the first instance, the argument runs
again from fine to coarse scale.
Let $f_{j+1}(x) = \Phi_{j+1}(x)\vec{s}_{j+1}$, then this function can be
projected onto the basis $\Phi_j(x)$. Even if power coefficients are projected
onto power coefficients, the projection is not unique, orthogonal
projection being just one of the possibilities. 
The reassembly of the refinement matrix $H_j$ from its factorisation, discussed
in this section, induces one particular projection.
Section \ref{subsec:Bsplinewavelets} will explain how to realise any other
projection using one further lifting step.
\end{revision}

A projection onto $\Phi_j(x)$ is characterised by a complementary basis
$\Psi_j(x)$, termed the wavelet basis, for which 
\begin{equation}
\Phi_{j+1}(x)\vec{s}_{j+1} = \Phi_j(x) \vec{s}_j + \Psi_j(x) \vec{d}_j.
\label{eq:FWTbasis1step}
\end{equation}
The expression (\ref{eq:FWTbasis1step}) is a basis transformation that
can be interpreted in two directions. From left to right, it represents the
projection, where the actual calculation of $\vec{s}_j$ and $\vec{d}_j$ still
has to be developed. In the other direction it
describes the reconstruction of the fine scale data from the coarse projection
$\Phi_j(x) \vec{s}_j$ and the residual $\Psi_j(x) \vec{d}_j$. 
The reconstruction includes the refinement given in the equation 
(\ref{eq:2scalematrix}). Indeed, take for the vectors
$\vec{s}_j$ and $\vec{d}_j$ a column of the matrices $I_{n_j}$ and
$O_{n_j}$, the identity and zero matrices of size $n_j \times n_j$. Then
the the vector $\vec{s}_{j+1}$ must be the corresponding column of $H_j$. In a
similar way one can take for $\vec{s}_j$ a zero column and for $\vec{d}_j$ a
canonical vector. The vector $\vec{s}_{j+1}$ is then the column of the matrix
$G_j$ in the wavelet equation
\begin{equation}
\Psi_j(x) = \Phi_{j+1}(x) G_j.
\label{eq:waveqmatrix}
\end{equation}
Substitution of the two-scale and wavelet equations (\ref{eq:2scalematrix})
and (\ref{eq:waveqmatrix}) into (\ref{eq:FWTbasis1step}) amounts to
\begin{equation}
\vec{s}_{j+1} = H_j \vec{s}_j + G_j \vec{d}_j.
\label{eq:IWT1step}
\end{equation}
The projection of $f_{j+1}(x)$ is found from the inverse of
(\ref{eq:IWT1step}). This inverse can be represented by the matrices
\begin{eqnarray}
\vec{s}_j & = & \widetilde{H}_j^T \vec{s}_{j+1},
\label{eq:FWT1stepHt}
\\
\vec{d}_j & = & \widetilde{G}_j^T \vec{s}_{j+1}.
\label{eq:FWT1stepGt}
\end{eqnarray}
The matrices $\widetilde{H}_j$ and $\widetilde{G}_j$ can be found from the
inversion of (\ref{eq:IWT1step}), which is formulated as the 
\emph{perfect reconstruction} property
\begin{equation}
H_j\widetilde{H}_j^T + G_j \widetilde{G}_j^T = I_{n_{j+1}}.
\label{eq:PR}
\end{equation}

The residual coefficients $\vec{d}_j$ are known as detail or wavelet
coefficients at scale $j$. The coarse scaling coefficients $\vec{s}_j$ can
further be processed into detail and scaling coefficients at scale $j-1$ and
so on. The multiscale transform from coefficients $\vec{s}_J$ at finest scale
$J$ into details at successive scales and scaling coefficients at a final,
coarse scale $\vec{s}_L$ is the forward wavelet transform or wavelet analysis.
It is carried out by repeated application of (\ref{eq:FWT1stepHt}) and
(\ref{eq:FWT1stepGt}). Likewise, (\ref{eq:IWT1step}) is one step in the inverse
wavelet transform or wavelet synthesis.

The forward transform matrix is denoted as $\widetilde{W}$. It maps the fine
scaling vector $\vec{s}_J$ onto the vector of coarse scaling coefficients
and multiscale details. Denoting the latter vector as $\vec{w}_L$ in the
definition 
\(
\vec{w}_L^T = \left[\begin{array}{cccc}
\vec{s}_L^T & \vec{d}_L^T & \ldots & \vec{d}_{J-1}^T
\end{array}\right],
\)
the forward wavelet transform is formalized as
\begin{equation}
\vec{w}_L = \widetilde{W} \vec{s}_J.
\label{eq:FWTmatrix}
\end{equation}
The inverse wavelet transform matrix is denoted by $W = \widetilde{W}^{-1}$.

The following theorem states that the lifting factorisation of a refinement
matrix $H_j$ can be used to find in a straightforward way a matrix $G_j^{[0]}$,
so that $H_j$ and $G_j^{[0]}$ constitute one step of an inverse wavelet
transform. Moreover, the forward transform matrices $\widetilde{H}_j^{[0]T}$ 
and $\widetilde{G}_j^T$ follow immediately as well, i.e., without any explicit
non-diagonal matrix inversion. 
The superscripts in $G_j^{[0]}$ and $\widetilde{H}_j^{[0]T}$ refer to the fact 
that these matrices follow in a natural way from the lifting factorisation of
$H_j$. Other matrices $G_j$ and $\widetilde{H}_j^T$ may appear in a perfect
reconstruction scheme (\ref{eq:PR}) with $H_j$ as well.
Theorem \ref{thm:finalU} will provide a lifting construction for all possible
$G_j$ and $\widetilde{H}_j^T$, given a refinement matrix $H_j$.
This construction does not affect $\widetilde{G}_j^T$, which is the reason for
not writing a superscript in that matrix.
\begin{thm} (wavelet transform from lifting factorisation)
If we can write $H_{j,e}^{[s]} = H_{j,e}^{[s+1]} - U_j^{[s+1]} H_{j,o}^{[s]}$,
as in (\ref{eq:primallifting}), with $H_{j,e}^{[s+1]} = D_j$ an invertible
matrix, then we can construct a wavelet transform containing $H_j^{[s]}$ as
refinement matrix.
In particular, we let $P_j^{[s+1]} = H_{j,o}^{[s]}D_j^{-1}$, and take as
forward transform
\begin{eqnarray}
\vec{s}_j & = & D_j^{-1}(\vec{s}_{j+1,e} + U_j^{[s+1]}\vec{s}_{j+1,o}),
\label{eq:finalfrwliftingstepU}
\\
\vec{d}_j & = & \vec{s}_{j+1,o} - P_j^{[s+1]}D_j\vec{s}_j.
\label{eq:finalfrwliftingstepP}
\end{eqnarray}
The synthesis or inverse transform consists of
\begin{eqnarray}
\vec{s}_{j+1,o} & = & \vec{d}_j + P_j^{[s+1]}D_j\vec{s}_j,
\label{eq:finalinvliftingstepP}
\\
\vec{s}_{j+1,e} & = & D_j\vec{s}_j - U_j^{[s+1]}\vec{s}_{j+1,o}.
\label{eq:finalinvliftingstepU}
\end{eqnarray}
\label{thm:lastliftingstep}
\end{thm}
\proof
From the inverse transform, we can check that $H_{j,o}^{[s]} = P_j^{[s+1]}D_j$,
while
$\vec{s}_{j+1,e} = D_j\vec{s}_j - U_j^{[s+1]}(\vec{d}_j + 
P_j^{[s+1]}D_j\vec{s}_j)$,
meaning that
$H_{j,e}^{[s]} = D_j - U_j^{[s+1]}P_j^{[s+1]}D_j = 
H_{j,e}^{[s+1]}- U_j^{[s+1]}H_{j,o}^{[s]}$.
\qed

The factorisation behind the transform is thus
\[
\left[\begin{array}{l}H_{j,e}^{[s]} \\ H_{j,o}^{[s]}\end{array}\right]
=
\left[\begin{array}{ll}I_{n_j} & -U_j^{[s+1]} \\ 0 & I_{n_j'}\end{array}\right]
\left[\begin{array}{ll}I_{n_j} & 0 \\ P_j^{[s+1]} & I_{n_j'}\end{array}\right]
\left[\begin{array}{l}D_j \\ 0\end{array}\right],
\]
where $n_j' = n_{j+1}-n_j$.

In principle, Theorem \ref{thm:lastliftingstep} is applicable to any
factorisation (\ref{eq:primallifting}). In practice, it becomes interesting 
as soon as $D_j$ is a diagonal matrix. Then the band structure for both forward
and inverse transform is under immediate control, as inverting the transform
requires no matrix inversion, except for the trivial case of a diagonal matrix.

Let $H_j$ be a general refinement matrix, so that $H_{j,e}$ has a larger
bandwidth than $H_{j,o}$, then a full factorisation is given by
\begin{equation}
H_j = \left[\begin{array}{l}H_{j,e} \\ H_{j,o}\end{array}\right]
=
\left(\prod_{s=1}^{u}
\left[\begin{array}{ll}I_{n_j} & -U_j^{[s]} \\ 0 & I_{n_j'}\end{array}\right]
\left[\begin{array}{ll}I_{n_j} & 0 \\ P_j^{[s]} & I_{n_j'}\end{array}\right]
\right)
\left[\begin{array}{l}D_j \\ 0\end{array}\right],
\end{equation}
where $u$ is the number of update steps.
The refinement matrix can be expanded into a full, invertible two-scale
transform by adding independent columns to the last factor, thus defining
a detail matrix 
\(
G_j^{[0]} = \left[\begin{array}{ll}
G_{j,e}^{[0]} & G_{j,o}^{[0]}\end{array}\right] 
\)
in the following factorisation
\begin{equation}
\left[\begin{array}{ll}
H_j & G_j^{[0]}
\end{array}\right]
=
\left[\begin{array}{ll}H_{j,e} & G_{j,e}^{[0]} \\ 
                       H_{j,o} & G_{j,o}^{[0]}
\end{array}\right]
=
\left(\prod_{s=1}^{u}
\left[\begin{array}{ll}I_{n_j} & -U_j^{[s]} \\ 0 & I_{n_j'}\end{array}\right]
\left[\begin{array}{ll}I_{n_j} & 0 \\ P_j^{[s]} & I_{n_j'}\end{array}\right]
\right)
\left[\begin{array}{ll}D_j & 0 \\ 0 & I_{n_j'}\end{array}\right].
\label{eq:IWTfactored}
\end{equation}
This is an inverse transform, i.e., the synthesis of fine scale coefficients.
The corresponding forward transform or analysis is denoted with the matrices
$\widetilde{H}_j^{[0]}$ and $G_j$ and its factorisation follows immediately
from the inversion of (\ref{eq:IWTfactored}), i.e.,
\begin{equation}
\left[\begin{array}{l}
\widetilde{H}_j^{[0]T}\\
\widetilde{G}_j^T
\end{array}\right]
=
\left[\begin{array}{ll}
H_j & G_j^{[0]}
\end{array}\right]^{-1}
=
\left[\begin{array}{ll}D_j^{-1} & 0 \\ 0 & I_{n_j'}\end{array}\right]
\left(\prod_{s=0}^{u-1}
\left[\begin{array}{ll}I_{n_j} & 0 \\ -P_j^{[q-s]} & I_{n_j'}\end{array}\right]
\left[\begin{array}{ll}I_{n_j} & U_j^{[q-s]} \\ 0 & I_{n_j'}\end{array}\right]
\right).
\label{eq:FWTfactored}
\end{equation}

For completeness, if $H_{j,o}$ has a larger bandwidth than $H_{j,e}$, then the
full factorisation of the refinement matrix is
\begin{equation}
H_j = \left[\begin{array}{l}H_{j,e} \\ H_{j,o}\end{array}\right]
=
\left[\begin{array}{ll}I_{n_j} & 0 \\ P_j^{[0]} & I_{n_j'}\end{array}\right]
\left(\prod_{s=1}^{u}
\left[\begin{array}{ll}I_{n_j} & -U_j^{[s]} \\ 0 & I_{n_j'}\end{array}\right]
\left[\begin{array}{ll}I_{n_j} & 0 \\ P_j^{[s]} & I_{n_j'}\end{array}\right]
\right)
\left[\begin{array}{l}D_j \\ 0\end{array}\right].
\end{equation}
The other matrices in the perfect reconstruction scheme (\ref{eq:PR}) follow
from this factorisation in a similar way as  in (\ref{eq:IWTfactored}) and
(\ref{eq:FWTfactored}).

In any case, the last step in the factorisation must be a prediction step.
Indeed, as follows from Theorem \ref{thm:lastliftingstep} and its
interpretation, the factorisation may stop as soon as $H_{j,e}^{[s+1]}$ is a
diagonal matrix. A diagonal matrix is not required for the last
$H_{j,o}^{[s+1]}$.

\begin{longversion}
Theorem \ref{thm:lastliftingstep} also implies to the following corollary
\citep{daubechies98:factorlifting}.
\begin{corollary}
A matrix $H_j$ can only qualify for a wavelet transform if all columns in
$H_{j,e}$ and $H_{j,o}$ are pairwise coprime.
\end{corollary}
Indeed, if the evens on a column form a vector which is a multiple of the odds,
then $D_j$ will contain a zero diagonal element in that column, making it a
singular matrix.
\end{longversion}

\subsection{Fast B-spline wavelet transforms}
\label{subsec:Bsplinewavelets}

So far, we have concentrated on finding the matrix $H_j$ in the two-scale
equation (\ref{eq:2scalematrix}) for a B-spline scaling basis $\Phi_j$.
The actual goal is, however, the design of a wavelet transform associated to
the B-spline basis. In particular, we want to choose an appropriate
wavelet basis $\Psi_j$ in (\ref{eq:FWTbasis1step}). From the wavelet equation
(\ref{eq:waveqmatrix}) it is clear that
properties of $\Psi_j(x)$ can be realised through an appropriate design of
$G_j$. The design of lifting steps (\ref{eq:IWTfactored})
in Section \ref{subsec:multiscaletransform} has
produced a matrix $G_j^{[0]}$ as a side effect, but it is unlikely that this
matrix realises the exact properties we have in mind. In particular, all
elements in $-U_j^{[s]}$ and $P_j^{[s]}$ are positive or zero, leading to the
conclusion that $G_j^{[0]}$ contains only non-negative entries. The functions
in $\Psi_j^{[0]}(x) = \Phi_{j+1}(x)G_j^{[0]}$, being linear combinations of
non-negative functions with a non-negative coefficients, cannot possibly have
zero integrals. Therefore, in a strict sense, these detail basis functions
cannot be termed wavelets.

The following theorem \citep{daubechies98:factorlifting} states that all
possible matrices $G_j$ can found one from the other by a single update
lifting step.
\begin{thm}(final update step)
Let $\vec{s}_{j+1} = H_j \vec{s}_j + G_j \vec{d}_j$ be a fine scale 
reconstruction with two-scale matrix $H_j$ and detail matrix $G_j$, and let
$\vec{s}_{j+1} = H_j \vec{s}_j^{[0]} + G_j^{[0]} \vec{d}$ be an alternative
scheme involving the same two-scale matrix and the same detail coefficients.
Both pairs $(H_j,G_j)$ and $(H_j,G_j^{[0]})$ belong to a perfect reconstruction
(\ref{eq:PR}) quadruple of matrices.
Then there exists un update operation $U_j$ so that
\begin{equation}
G_j = G_j^{[0]} - H_j U_j.
\label{eq:finalU}
\end{equation}
As a consequence, the updated scaling coefficients $\vec{s}_j$ can be
found from $\vec{s}_j^{[0]}$ and $\vec{d}_j$:
\begin{equation}
\vec{s}_j^{[1]} = \vec{s}_j^{[0]} + U_j \vec{d}_j.
\label{eq:finalupdatestep}
\end{equation}
\label{thm:finalU}
\end{thm}
\proof
The proof is straightforward by the construction
\begin{equation}
U_j = \widetilde{H}_j^TG_j^{[0]}.
\label{eq:Ufinal}
\end{equation}
Perfect reconstruction and Definition (\ref{eq:Ufinal}), amount to
\[
\left[\begin{array}{l}
\widetilde{H}_j^T\\
\widetilde{G}_j^T
\end{array}\right]
\left[\begin{array}{ll}
H_j & G_j^{[0]}
\end{array}\right]
=
\left[\begin{array}{ll}
I_{j+1,e} & U_j\\0 & I_{j+1,o}
\end{array}\right],
\mbox{ and to }
\left[\begin{array}{l}
\widetilde{H}_j^T\\
\widetilde{G}_j^T
\end{array}\right]
\left[\begin{array}{ll}
H_j & G_j
\end{array}\right]
=
\left[\begin{array}{ll}
I_{j+1,e} & 0\\0 & I_{j+1,o}
\end{array}\right].
\]
These two expressions can be combined into
\[
\left[\begin{array}{l}
\widetilde{H}_j^T\\
\widetilde{G}_j^T
\end{array}\right]
\left[\begin{array}{ll}
H_j & G_j^{[0]}
\end{array}\right]
-
\left[\begin{array}{ll}
I_{j+1,e} & U_j\\0 & I_{j+1,o}
\end{array}\right]
\left[\begin{array}{l}
\widetilde{H}_j^T\\
\widetilde{G}_j^T
\end{array}\right]
\left[\begin{array}{ll}
H_j & G_j
\end{array}\right]
=
\left[\begin{array}{ll}
0 & 0 \\ 0 & 0
\end{array}\right],
\]
which includes (\ref{eq:finalU}).
\qed

In the factored wavelet transforms of (\ref{eq:IWTfactored}) and
(\ref{eq:FWTfactored}) the operation $U_j$ occurs at the end of the forward
transform. Formally, this adds
\begin{equation}
\left[\begin{array}{ll}H_{j,e} & G_{j,e} \\ 
                       H_{j,o} & G_{j,o}
 \end{array}\right]
=
\left[\begin{array}{ll}H_{j,e} & G_{j,e}^{[0]} \\ 
                       H_{j,o} & G_{j,o}^{[0]}
 \end{array}\right]
\left[\begin{array}{ll}I_{n_j} & -U_j \\ 0 & I_{n_j'}\end{array}\right]
\end{equation}
to (\ref{eq:IWTfactored}) and, obviously,
\begin{equation}
\left[\begin{array}{ll}H_{j,e} & G_{j,e} \\ 
                       H_{j,o} & G_{j,o}
 \end{array}\right]^{-1}
=
\left[\begin{array}{ll}I_{n_j} & U_j \\ 0 & I_{n_j'}\end{array}\right]
\left[\begin{array}{ll}H_{j,e} & G_{j,e}^{[0]} \\ 
                       H_{j,o} & G_{j,o}^{[0]}
 \end{array}\right]^{-1}
\end{equation}
to (\ref{eq:FWTfactored})

Unlike the lifting steps in Section \ref{subsec:multiscaletransform}, the
matrix $U_j$ can have more than two non-zeros in each of its columns.

Given a matrix $G_j^{[0]}$, for instance from the factoring in Sections
\ref{subsec:multiscaletransform} and \ref{subsec:designliftingsteps}, the
design of the final wavelet matrix $G_j$ proceeds through the design of the
final update $U_j$. A typical design property is to impose vanishing moments,
that is $\int_{-\infty}^\infty \psi_{j,k}(x) x^m = 0$ for all $j$ and $k$ and
for $m=0,\ldots,p$.
Define the moment vectors of $\Phi_j(x)$ and $\Psi_j(x)$ as
\begin{eqnarray}
M_{j;m} & = & \int_{-\infty}^\infty \Phi_j^T(x) x^m dx,
\label{def:momentsPhij}
\\
O_{j,m} & = & \int_{-\infty}^\infty \Psi_j^T(x) x^m dx,
\label{def:momentsPsij}
\end{eqnarray}
and let 
\(
M_j^{(p)} = \left[
\begin{array}{cccc}M_{j;0} & M_{j;1} & \ldots & M_{j;p-1}\end{array}
\right],
\)
and similarly for $O_j^{(p)}$.
As $\Psi_j^{[0]}(x) = \Phi_{j+1}(x)G_j^{[0]}$, the preliminary moments can be
computed throughout the lifting scheme culminating into the expression
$O_j^{(p)[0]} = G_j^{[0]T} M_{j+1}^{(p)}$. Similarly we find $M_j^{(p)} = H_j^T
M_{j+1}^{(p)}$. As the final update defines the basis function $\Psi_j(x) =
\Phi_{j+1}(x)G_j = \Psi_j^{[0]}(x) - \Phi_j(x) U_j$, we have $O_j =
O_j^{(p)[0]} - U_j^T M_j^{(p)}$. Imposing that $O_j = [0]$ amounts to the
equation
\begin{equation}
O_j^{(p)[0]} = U_j^T M_j^{(p)}.
\label{eq:Uvanmom}
\end{equation}

For a wavelet in the strict sense of the definition, we need to impose
the vanishing moment condition for $p=1$. Higher vanishing moments 
are often not so useful in statistical applications, and if they turn out to be
beneficial, then often it is a better idea to optimise for the wanted benefits
in a more explicit way. In particular, for use in statistics, it is better to
explicitly impose that the transform is as close as possible to being
orthogonal.
All this is further developed in Section \ref{sec:estBsplinewavelets}.

\subsection{The non-decimated B-spline wavelet transform}
\label{subsec:RWT}

The non-decimated wavelet transform is a redundant data decomposition that has
$n$ wavelet coefficients at each scale, where $n$ is the size of the input
vector. Moreover, if the elements of the input vector is shifted, then the
coefficients at each scale are shifted in the same way. This is the translation
invariance property. Translation invariance is impossible in a decimated
transform. Indeed, the decimation takes place in the the even-odd partitioning.
As evens and odds play different roles in the subsequent lifting steps, a shift
in the input vector leads to a different role for each input element, and thus
a different outcome. The fact that shifted inputs lead to outcomes with
different values may, for obvious reasons, complicate the analysis or
processing of the wavelet coefficients.

Each step in the non-decimated transform starts from $n$ scaling coefficients
at scale $j+1$ and produces $n$ scaling coefficients at scale $j$ together with
$n$ wavelet coefficients. At the finest scale, i.e., for $j=J-1$, the decimated
scaling coefficients fill up $\ceil{n/2}$ values of the non-decimated
expansion, while $\floor{n/2}$ values of the non-decimated wavelet coefficients
at scale $J-1$ come from the decimated transform. In order to complete the
other half of the coefficients, the same transform is carried out switching the
roles of evens and odds. This can be realised by defining the shifted knots
$x_{J,k}^{[1]} = x_{J,k-1}$, while we set $x_{J,k}^{[0]} = x_{J,k} = x_k$ for
the original vector of knots. After the first step, the shifted vector of knots
has generated an alternative decimated set of knots, which contains the evens
of the shifted vector, that is, the odds of the original vector. We denote the
new vector at scale $J-1$ by
$x_{J-1,k}^{[2]} = x_{J,2k}^{[1]} = x_{J,2k+1} = x_{2k+1}$, 
while the original decimated vector is now denoted as
$x_{J-1,k}^{[0]} = x_{J,2k}^{[0]} = x_{J,2k} = x_{2k}$.
Both vectors can be shifted for use in the second step, thereby defining two
more vectors
$x_{J-1,k}^{[1]} = x_{J-1,k-1}^{[0]} = x_{J,2k-2}^{[0]}$ and
$x_{J-1,k}^{[3]} = x_{J-1,k-1}^{[2]} = x_{J,2k-2}^{[1]}$.
All four vectors are used in the second step, leading us to scale $J-2$.
In general, the non-decimated transform at scale $j$ consists of $2^{J-j}$
decimated transforms, each defined by a vector of knots
$x_{j,k}^{[2a+b]} = x_{j+1,2k-2b}^{[a]} =
x_{(k-b)2^{J-j}+\mathrm{rem}(a,2^{J-j+1})}$.

\section{The experimental approximation error of a spline wavelet decomposition}
\label{sec:apprxerr}

The following subsections discuss the accuracy of approximation in a B-spline
wavelet basis on irregular knots. The importance of the approximation error in
statistical applications lies in the fact that approximation error is a
source of estimation bias. From the discussion below, it turns out that
the common practice in wavelet based noise reduction to take observations as
fine scaling coefficients can not be adopted when the observations are not
equidistant.

\subsection{Linear approximation error in a homogeneous dyadic refinement}
\label{subsec:apprxerrlinear}

Let $f(x)$ be a smooth function on $[0,1]$, 
i.e., a function with at least $\widetilde{p}$
continuous and bounded derivatives, or a function which is uniformly
Lipschitz-$\alpha$ continuous, with $\alpha\geq \widetilde{p}$, as defined in
\citep[Definition 6.1]{mallat01:wavelettour2nd}. This function can be
approximated by a linear combination $f_J(x)$ of the B-spline scaling functions
at level $J$ defined on the equidistant knots $x_{J,k} = k2^{-J}$, leading to
an expansion as in (\ref{eq:expansionBsplinebasis}), substituting $j$ by $J$.
As in Section \ref{subsec:multiscaletransform}, index $J$ refers to the finest
scale, this is the scale at which the observations take place.
It can be proven that both the pointwise approximation error $f(x)-f_J(x)$ 
\citep[(3.5) and (3.9)]{sweldens94:quadrature} and the mean squared error
$\norm{f(x)-f_J(x)}$ \citep[Theorem 7.5, page 230]{strang96:filterbanks}
are of the order $\mathcal{O}\left(2^{-J\widetilde{p}}\right)$. The
approximation $f_J(x)$ is not unique. It can, for instance, be defined in
accordance to a subsequent wavelet decomposition that is applied to the
approximation. In principle, the wavelet transform of the approximation
does not interact with the construction of the approximation at finest scale.
Nevertheless, each wavelet system $\psi_{j,k}(x)$ fixes a dual scaling basis
$\widetilde{\varphi}_{J,k}(x)$ which can be used to define a projection 
\(
f_J(x) =
\sum_{k=0}^{2^J-1}\dotp{f(x),\widetilde{\varphi}_{J,k}(x)}\varphi_{J,k}(x).
\)
Alternatively, if $\varphi_{J,k}(x)$ does not form an orthogonal basis, then
$f_J(x)$ can also be the orthogonal projection onto the basis
$\varphi_{J,k}(x)$. Suggesting yet another possibility, similar approximation
results are also available for interpolating splines on regular point sets
\citep{dubeau95:errorboundsspline}. In all these cases, we
conclude that a resolution of $2^{-J}$ is enough to obtain an accuracy of
$2^{-J\widetilde{p}}$, thanks to the smoothness of the function $f(x)$. To the
best of my knowledge, nearly no such results exist for irregular knots. One of
the difficulties is that little is known about the dual functions
$\widetilde{\varphi}_{J,k}(x)$. On regular knots, all these functions are
translations and dilation of a single dual father function, allowing us in a
fairly easy way to establish a general upper bound for the approximation error.

The accuracy of order $2^{-J\widetilde{p}}$ for resolution $2^{-J}$ does not
hold for the approximation obtained by taking function values as scaling
coefficients, i.e.,
\begin{equation}
\overline{f}_J(x) = \sum_{k=0}^{2^J-1} f(x_{J,k}) \varphi_{J,k}(x).
\label{eq:functionvaluesscalingcoefs}
\end{equation}
Instead, it can be proven \citep[Theorem 2.4]{sweldens94:quadrature} that
$\overline{f}_J(x)$ is within an $L_2$-distance 
$\mathcal{O}\left(2^{-J\widetilde{p}}\right)$ from the blurred function
\begin{equation}
\barbar{f}_J(x) = \sum_{k=0}^{2^J-1} \varphi_{J,0}(x_{J,k})f(x-x_{J,k}),
\label{eq:waveletcrime}
\end{equation}
whose approximation error, it its turn, is of the order
$\norm{\barbar{f}_J(x)-f(x)}_2 = \mathcal{O}(2^{-J})$, at least if the
$\varphi_{J,k}(x)$ are normalised so that
\(
\sum_{k=0}^{2^J-1} \varphi_{J,k}(x) = 1.
\)
The blurring effect thus neutralises all
benefits from the $\widetilde{p}$ vanishing moments for linear approximation of
smooth function by a refinable basis. This phenomenon is known as the
\emph{wavelet crime} \citep{strang96:filterbanks}. If, however, the scaling
functions are interpolating in the sense that $\varphi_{0,0}(k) = 0$ unless
$k=0$, then the approximation accuracy is restored. This is the case for the
wavelets that follow from the Deslauries-Dubuc refinement scheme
\citep{deslauriers87:interpolation,deslauriers89:symmetric,donoho99:deslauries,sweldens96:athome},
and the advantage is preserved if the Deslauries-Dubuc refinement takes place
on a non-equidistant point set. Although Deslauries-Dubuc refinement may suffer
in other aspects from non-equidistance, the immediate use of function values at
the input is an important benefit for this scheme.
As an alternative for interpolating scaling functions, one can
impose that the projection coefficients are close to the function values, i.e.,
\(
\dotp{f(x),\widetilde{\varphi}_{J,k}(x)} = f(x_{J,k}) +
\mathcal{O}\left(2^{-J\widetilde{p}}\right).
\)
This is realised by imposing that $\dotp{x^q,\widetilde{\varphi}_{J,k}(x)} = 0$
for $q \in \{1,\ldots,\widetilde{p}-1\}$.
If the basis is orthogonal, then $\widetilde{\varphi}_{J,k}(x) =
\varphi_{J,k}(x)$, and the development of the condition leads to the class of
coiflets \citep{daubechies93:coiflets}.

\subsection{Nonlinear approximation, compression and thresholding}
\label{subsec:apprxerrnonlinear}

When $f(x)$ contains jumps, cusps, or other singularities, any approximation
as in (\ref{eq:expansionBsplinebasis}) may have a local error of order
$\mathcal{O}(1)$ near the singularities. More precisely, if the interval
$[x_{J,k},x_{J,k+1}]$ contains a singularity, then for any point $x \in
[x_{J,k},x_{J,k+1}]$, the pointwise approximation error is
$|f_J(x)-f(x)| = \mathcal{O}(1)$. As the location of the singularity is
only known up to the resolution of the observation, the local error contributes
to the total $L_2$-error at a rate equal to the resolution of the observation,
no matter how accurately the smooth parts in between are approximated.
Therefore, when $f(x)$ is piecewise smooth, then the
resolution of observations, $J$, is often taken much finer than necessary for
the application. Next, a wavelet decomposition is applied and all fine scale
detail coefficients up to a level $L$ are omitted, except for those that
correspond to the singularities (typically the large ones). Taking
$J>L\widetilde{p}$ ensures that the error in catching the singularities does
not exceed the error from the smooth part approximation.

In this nonlinear thresholding scheme, the error from using function values as
scaling coefficients at scale $2^{-J}$ has to be compared to the smooth
approximation error of $2^{-L\widetilde{p}}$. Setting $J>L\widetilde{p}$ thus
also ensures that function values as finest scaling coefficients poses no
problem, \emph{at least if} the thresholding in between causes no additional
error. In particular, if $f(x)$ is a polynomial, then all wavelet coefficients
of a proper approximation $f_J(x)$ are zero. In that case, the detail
coefficients of $\overline{f}_J(x)$ are also zero, because of the following
result \citep[Theorem 7.4, p.241]{mallat01:wavelettour2nd}.
\begin{thm}
Let $\varphi(x)$ be a father scaling function. Then for all $q \in
\{0,1,\ldots,\widetilde{p}-1\},$ the function
\[
\overline{v}^{\left[\widetilde{p},q\right]}(x) =
\sum_{k \in \ZZ} k^q \varphi(x-k),
\]
is a polynomial of degree $q$ if and only if for all $q \in
\{0,1,\ldots,\widetilde{p}-1\},$ there exist a sequence of
coefficients $\widetilde{x}_k^{[\widetilde{p},q]}$ so that
\[
x^q = \sum_{k \in \ZZ} \widetilde{x}_k^{\left[\widetilde{p},q\right]}
\varphi(x-k).
\]
\label{thm:StrangFix}
\end{thm}
This result has the following interpretation: if the basis consists of
translations of a single father function along equispaced knots, and if the
basis reproduces polynomials up to degree $\widetilde{p}-1$, then the
coefficients for the decomposition of that polynomial can be found as the
evaluation of another polynomial in the knots. Conversely, if the function
values of any polynomial used as coefficients lead to a polynomial, then all
polynomials can be reproduced within this basis.

Replacing the scaling coefficients of a polynomial by the function values thus
defines a new polynomial, whose detail coefficients are also zero, at least if
the scaling functions have $\widetilde{p}$ vanishing moments to reproduce
polynomials up to degree $\widetilde{p}-1$. Thresholding up to level $L$
introduces then no additional error. This motivates the
$\widetilde{p}$ vanishing moments in compression and denoising, even if at the
finest scale, they are not exploited when function values are plugged in as
scaling coefficients. As a conclusion, on an equidistant set of knots, and in a
nonlinear wavelet method, the wavelet crime can be forgiven.

\subsection{Nonlinear approximation on non-equidistant knots}
\label{subsec:apprxerrirr}

In an approximation using non-interpolating scaling functions, such as
B-splines, on a non-equidistant set of knots, the wavelet crime is
unforgivable. It has two unpleasant effects that do not occur in the equispaced
case. First, the approximation error may propagate towards the coarser levels.
Indeed, Theorem \ref{thm:StrangFix} is not applicable. From
(\ref{eq:linearcoefsBspline}), it can be seen that the coefficients
representing the identity function on the grid of knots cannot be
retrieved as the evaluation of a polynomial in the knots. Conversely, when we
take the knot values $x_{J,k}$ as scaling coefficients, they will not be
recognised as coming from a smooth function. We can write
\(
x_{J,k} =
\widetilde{x}_{J,k}^{\left[\widetilde{p},1\right]} + \varepsilon_{J,k},
\)
where the error term $\varepsilon_{J,k}$ depends on the configuration of the
knots. 
Thresholding the detail coefficients that follow from such errors
$\varepsilon_{J,k}$ may have an error reducing effect. Depending on the 
positions of the knots, it may also lead to an increase of the
approximation error up to the order of $\mathcal{O}(\Delta_{L})$,
thereby undoing all benefits from $\varphi_{j,k}(x)$ having $\widetilde{p}$
vanishing moments.
The second unpleasant effect of taking function values as scaling coefficients
is a visual one: it leads to a decrease in smoothness of the reconstruction at
the fine level.
Indeed, the approximation $\overline{f}_J(x)$, reflects the irregular spacing
of the knots $x_{J,k}$, in a way similar to a decomposition and reconstruction
that would ignore the non-equidistant locations, assuming that
$x_{J,k} = k/2^J$.

Therefore, before carrying out the actual wavelet analysis, we need to
find a vector of scaling coefficients $\vec{s}_J$ for which $f(x) =
\Phi_J^{\left[\widetilde{p}\right]}(x)\vec{s}_J$.
As we are given the function values $f(x_{J,k})$ for
$k=0,\ldots,n-1$, we have to solve the set of equations $f(x_{J,k}) =
\sum_{l=0}^{n-1}s_{J,l}\varphi_{J,l}^{\left[\widetilde{p}\right]}(x_{J,k})$.
In matrix notation, this is
\(
\vec{f}_J = \bm{\Phi}_J \vec{s}_J,
\)
where $\vec{f}_J = [\ldots,f(x_{J,k}),\ldots]$ is the vector of observations
and where the matrix $\bm{\Phi}_J$ has entries 
$\bm{\Phi}_{J;k,l} = \varphi_{J,l}^{\left[\widetilde{p}\right]}(x_{J,k})$. For
the sake of readability, the superscripts are omitted in the notation of the
matrix.
The evaluations $\varphi_{J,l}^{\left[\widetilde{p}\right]}(x_{J,k})$ are
carried out using the recursion (\ref{eq:Bsplinerecursion}) in Definition
\ref{def:Bspline}. All B-splines are zero in the first and last knots, i.e.,
$\varphi_{J,l}^{\left[\widetilde{p}\right]}(x_{J,k}) = 0$ for $k \in
\{0,n-1\}$, and for all $l \in \{0,\ldots,n-1\}$. So, the first and last row in
$\bm{\Phi}_J$ contain zeros only, and no vector $\vec{s}_J$ can reconstruct
the values $f(x_{J,0})$ and $f(x_{J,n-1})$.
In order to be able to reconstruct $f(x_{J,0})$ and
$f(x_{J,n-1})$, we add two artificial knots left and right from the interval
$[x_{J,0},x_{J,n-1}]$. 
By adding the two columns for the additional corresponding
B-splines, but not the zero rows for these new end knots, we get an 
$n \times (n+2)$ matrix $\overline{\bm{\Phi}}_J$
and a vector $\overline{\vec{s}}_J$ of length $(n+2)$.
The system $\overline{\bm{\Phi}}_J\overline{\vec{s}}_J = \vec{f}_J$ being 
indeterminate system, we look for a solution that behaves well if $\vec{f}_J$
is observed with noise. For reasons of superposition, the transformation
$\overline{\vec{s}}_J = S_J \vec{f}_J$ should be linear.
There is no need for a noise reducing effect in the transformation, since we 
want to keep the noise reduction for subsequent wavelet analysis, which is
much better equipped for the task, especially when $\vec{f}_J$ has
singularities. On the other hand, independent, homoscedastic noise should stay
more or less homoscedastic. Therefore, the matrix $S_J$ should be as close as
possible to the identity matrix. In particular, its singular values should be
close to 1. This would motivate to take $S_J = \overline{\bm{\Phi}}_J^T
(\overline{\bm{\Phi}}_J\overline{\bm{\Phi}}_J^T)^{-1}$. Simulation studies
show, however, that even in this optimal case, the singular values are much too
large for use in practice. Moreover, the outcome $\overline{\vec{s}}_J$ would
be an exact solution of the indeterminate system, but not necessarily the
solution that the subsequent wavelet analysis expects. In particular, suppose
that $\vec{f}_J$ comes from a polynomial of degree less that $\widetilde{p}$,
then the wavelet detail coefficients should be zero. This happens only if the
operation $S_J$ would deliver the right power coefficients, as given in Theorem
\ref{thm:x^qcoefBspline}. There is no guarantee that this would
happen. For these two reasons, no exact solution of the system
$\overline{\bm{\Phi}}_J\overline{\vec{s}}_J = \vec{f}_J$ is satisfactory for
use in practice.

An approximative solution can be found using some sort of regularisation.
Tikhonov regularisation would minimise
$\norm{\overline{\bm{\Phi}}_J\overline{\vec{s}}_J - \vec{f}_J}_2^2 + \alpha
\norm{\overline{\vec{s}}_J}_2^2$, for some appropriate value of $\alpha$. As
the objective is to control but not to reduce the noise variance, the value of
$\alpha$ would be smaller than in ridge regression. While the regularisation
controls the variance effectively, there is no control at all on the
propagation of approximation error, i.e., the bias, through the subsequent
wavelet decomposition. The same remark holds for other general regularisation
methods, such as the Landweber iteration scheme for ill posed linear inverse
problems \citep{landweber51:iterative,daubechies04:iterativethresholding}.

A good compromise for $S_J$ should not be too restrictive, i.e., it should not
impose the perfect reconstruction $\overline{\bm{\Phi}}_JS_J\vec{f}_J =
\vec{f}_J$ for any vector $\vec{f}_J$. Instead, we focus on the power
functions. Let $\bm{X}^{[\widetilde{p}]}$ denote the matrix with entries
$(\bm{X}^{[\widetilde{p}]})_{k,q+1} = x_{J,k}^q$ for
$q=0,\ldots,\widetilde{p}-1$ and let
$\widetilde{\bm{X}}^{[\widetilde{p}]}$ the matrix with the corresponding
coefficients in the fine scale B-spline basis, i.e.,
$(\widetilde{\bm{X}}^{[\widetilde{p}]})_{k,q} =
\widetilde{x}_{J,k}^{\left[\widetilde{p},q\right]}.$ 
Then imposing
\begin{equation}
S_J \bm{X}^{[\widetilde{p}]} = \widetilde{\bm{X}}^{[\widetilde{p}]},
\label{eq:S_Jconstraint}
\end{equation}
ensures that all polynomials of degree less than $\widetilde{p}$ are
represented exactly and with the coefficients that are recognised by the
subsequent wavelet transform as coming from a polynomial.
Expression (\ref{eq:S_Jconstraint}) formulates $\widetilde{p}$ conditions for
each row $k$ of $S_J$. Let $\partial k$ denote the set of indices $l$ for
which we choose $S_{J;k,l} \neq 0$, then taking $\# \partial k \geq
\widetilde{p}+1$, allows us to satisfy these constraints, leaving one or a few
free parameters to control the variance of the scaling coefficients. In
particular, an objective can be to take $S_J$ as close as possible to the
identical transform, extended with two zero columns for the two additional
knots at the end points of the interval. We minimise the Frobenius norm
$\norm{S_J - \overline{I}_J}_F$, where
\(
\overline{I}_J = 
\left[\begin{array}{ccc}
\vec{0}_n & I_{n_J} & \vec{0}_n
\end{array}\right].
\)
On the level of the $k$th row, this amounts to the constraint optimisation
problem
\begin{equation}
\min_{S_{J;k,\partial k}} 
\norm{S_{J;k,\partial k} - \delta_{k+1,\partial k}}_2^2,
\end{equation}
subject to
\begin{equation}
\sum_{l \in \partial k} S_{J,k,l} x_{J,l}^q =
\widetilde{x}_{j,k}^{\left[\widetilde{p},q\right]},
\end{equation}
for $q = 0,\ldots,\widetilde{p}-1$.

\section{Estimation in a B-spline wavelet basis}
\label{sec:estBsplinewavelets}

The spline wavelets proposed by 
\citet*{coh-dau-fea92:biorthogonal}
devote all free parameters in $G_j$
to vanishing moments $p$. This corresponds to a lifting scheme where the final
update step is found by the system in (\ref{eq:Uvanmom}). Imposing one primal
vanishing moment, i.e., (\ref{eq:Uvanmom}) with $p=1$, is necessary to have
wavelets in the strict sense of the word, basis functions that fluctuate around
zero so that their integral is zero. Basis functions with zero integral are
indispensable for the representation of square-integrable functions. Indeed,
when the basis functions have nonzero integrals, there exist nontrivial
approximations of the zero function that converge in quadratic norm
\citep[page 93]{jansen05:book}. Basis functions with nonzero integrals are
useful on subspaces of the square-integrable functions, defined by additional
smoothness conditions. These conditions exclude functions with jumps, which are
typically the functions of interest in a wavelet analysis.

\begin{figure}[t!]
\begin{tabular}{cc}
\includegraphics[width=0.49\textwidth]{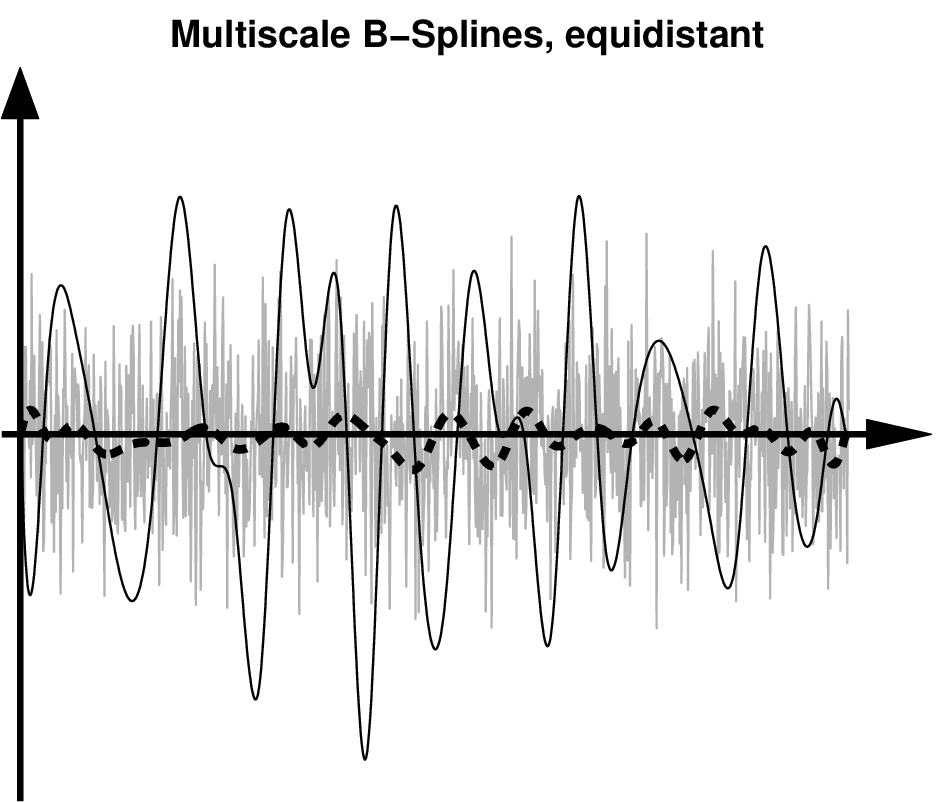} &
\includegraphics[width=0.49\textwidth]{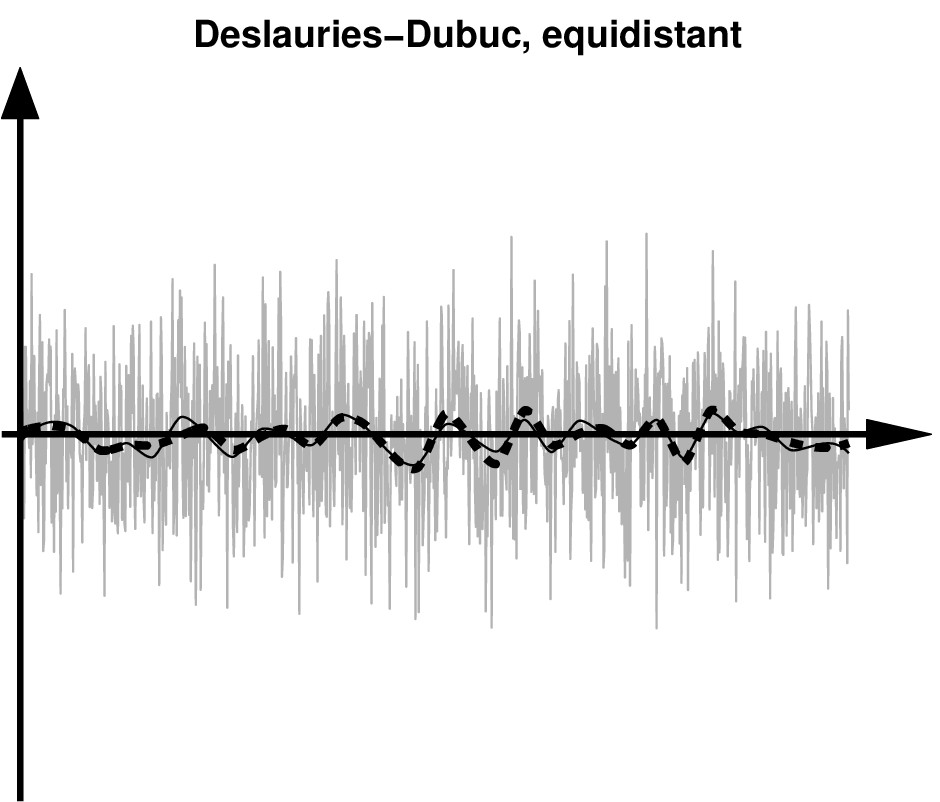} \\
\includegraphics[width=0.49\textwidth]{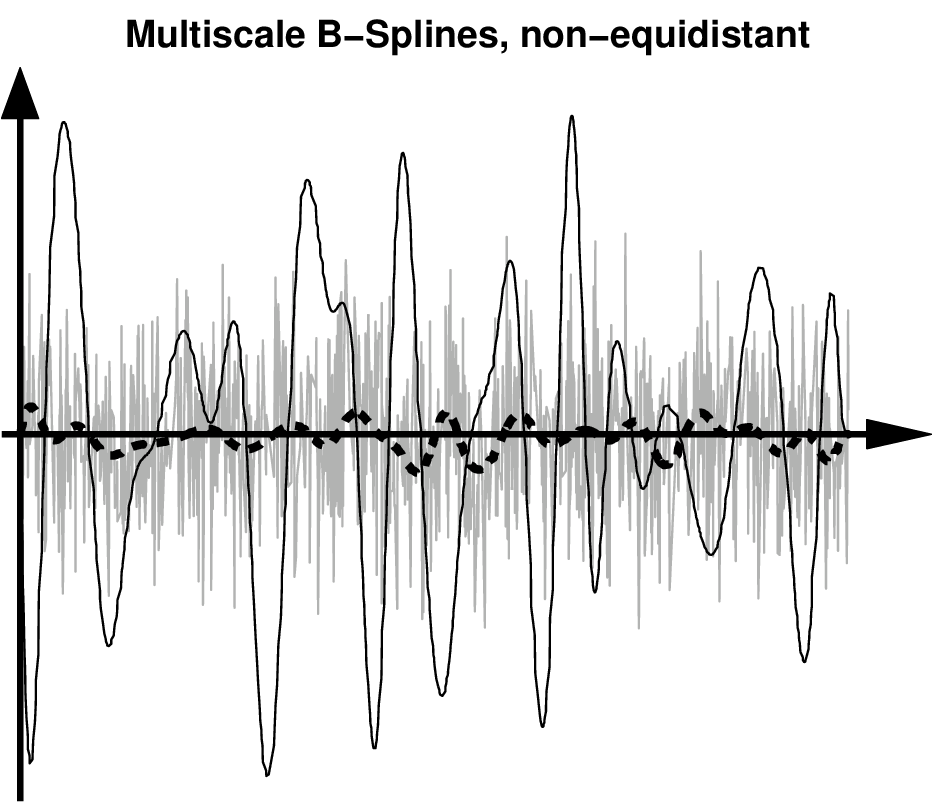} &
\includegraphics[width=0.49\textwidth]{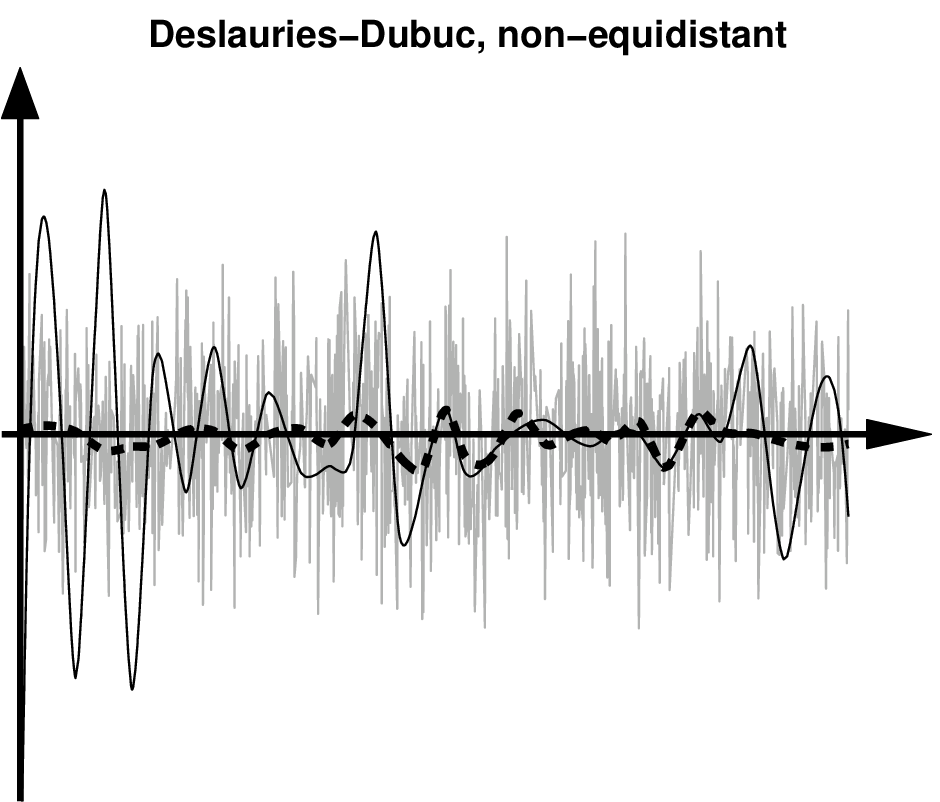}
\end{tabular}
\caption{Noise reduction using linear diagonal selection in a multiscale cubic
B-spline basis and in a Deslauries-Dubuc scheme with multiscale cubic
interpolation: reconstruction with all detail coefficients replaced by zero.
Equidistant and non-equidistant observations.
In grey line the observations, in black solid line the reconstructions from
projections with two primal vanishing moments, in thick black dashed line the
reconstructions from the orthogonal projections onto the B-spline or
Deslauries-Dubuc bases.}
\label{fig:expvar}
\end{figure}

On the other hand, the experiment in Figure \ref{fig:expvar} illustrates that 
using all free parameters for a maximum number of primal vanishing moments may
not be the best choice in a context of function estimation. This holds in
particular on non-equispaced data, irrespective of the wavelet family. For
spline wavelets, it holds also for data observed in equidistant points. The
experiment suggests that specific design criteria are necessary for wavelets in
statistical applications.

The experiment starts from $n_J = 1000$
fine scaling coefficients that are independently and
identically distributed random variables with zero mean $\vec{\varepsilon}_J$.
Using $\widetilde{W}$ from (\ref{eq:FWTmatrix}), define the transformed random
vector $\vec{\eta}_L = \widetilde{W}\vec{\varepsilon}_J$.
The wavelet coefficient vector $\vec{\eta}_L$ consists of coarse scaling
coefficients plus detail coefficients at each scale,
\(
\vec{\eta}_L^T = \left[\begin{array}{cccc}
\vec{\varepsilon}_L^T & \vec{\delta}_L^T & \ldots & \vec{\delta}_{J-1}^T
\end{array}\right].
\)
Let $D_L$ be the linear diagonal selection operation that keeps the coarse
scaling coefficients and discards all detail coefficients, i.e., 
\(
D_L \vec{\eta}_L = \left[\begin{array}{cccc}
\vec{\varepsilon}_L^T & \vec{0}_L^T & \ldots & \vec{0}_{J-1}^T,
\end{array}\right]^T
\)
and consider the reconstruction
$WD_L \vec{\eta}_L^J = (WD_L\widetilde{W}) \vec{\varepsilon}_J = 
(WD_LD_L\widetilde{W}) \vec{\varepsilon}_J =
(W_L\widetilde{W}_L)\vec{\varepsilon}_J$.
Here $\widetilde{W}_L = D_L\widetilde{W}$ and $W_L = WD_L$.
Note that
\begin{equation}
W_L = \prod_{i=1}^{J-L} H_{J-i} = H_{J-1}H_{J-2}\ldots H_{L+1}H_L,
\label{eq:WL}
\end{equation}
while
\begin{equation}
\widetilde{W}_L = \prod_{j=L}^{J-1} \widetilde{H}_j^T.
\label{eq:WtL}
\end{equation}

Since all its steps are linear operations, the reconstruction 
$(W_L\widetilde{W}_L)\vec{\varepsilon}_J$ is unbiased. All
observed errors are due to the variance of the reconstruction.
Each plot in Figure \ref{fig:expvar} compares a non-orthogonal projection
$\mathrm{P}_L = W_L\widetilde{W}_L$ with the orthogonal 
$\mathrm{P}_{L\perp} = W_L(W_L^TW_L)^{-1}W_L^T$. 
The projection matrices $\mathrm{P}_L$ should not be confused with the
prediction matrices $P_j^{[s]}$ in Sections \ref{subsec:Bspline2scales},
\ref{subsec:multiscaletransform}, and \ref{subsec:designliftingsteps}.
The non-orthogonal projection is constructed level by level, where in each
level the lifting factorisation of the cubic B-spline basis is completed with a
final update matrix $U_j$ that has one nonzero off-diagonal next to a nonzero
diagonal. All nonzero elements in $U_j$ are filled in by imposing two primal
vanishing moments, i.e., $O_{j,m} = \vec{0}_j$ for intermediate scales and for
$m\in\{0,1\}$, where $O_{j,m}$ is defined by (\ref{def:momentsPsij}).

The reconstruction from the Deslauries-Dubuc projection with two vanishing 
moments shows a large variance on the non-equidistant data, but not on the
equidistant equivalent. Reconstruction from a projection with two vanishing
moments onto a B-spline basis shows a large variance in both the equidistant
and non-equidistant cases. Moreover, further experiments reveal that the
variance increases when the B-spline wavelet transform involves more scales.
In the Deslauries-Dubuc scheme, there appears to be no such multiscale
deterioration. Large variances in that scheme are rather isolated, although
possibly devastating, effects from local irregularities in the non-equidistant
grid of observations.

Although the large variances are clearly an effect of the non-orthogonality of
the transform, the classical numerical condition number sheds no light on the
problem.
Indeed, on an equispaced, dyadic grid of $n_J=2048$ knots, the condition
number of the forward Cohen-Daubechies-Feauveau cubic B-spline wavelet with two
primal vanishing moments turned out to be $71.9$, and this number is fairly
independent from $n_J$. On the other hand the condition number
of a cubic Deslauries-Dubuc scheme with also two primal vanishing moments is
dependent on $n_J$ and for $n_J=2048$ it equals $203.7$. Nevertheless, the
Deslauries-Dubuc scheme on equispaced data shows no problems with large
variances. The numerical condition of the wavelet transform is therefore not an
adequate description of the non-orthogonality of the transform for
statistical applications. The same conclusion holds for the Riesz constants 
\cite{mallat01:wavelettour2nd} in a biorthogonal transformation.

As in B-spline wavelet transforms the increasing variances across scales occur
also on equidistant data sets, it is possible to describe them using Fourier
transforms. Given a vector after projection
$\vec{\varepsilon}_L = \mathrm{P}_L\vec{\varepsilon}$, define its Fourier
transform
$Y(\omega) = {1 \over n} \sum_{\ell=1}^n \varepsilon_{L,i} e^{-\ell i \omega}$.
By carefully checking the effect of sub- and up-sampling on the Fourier
transform, it is fairly straightforward to prove that for $\varepsilon$
uncorrelated and homoscedastic, it holds that
\begin{equation}
E\left|Y(\omega)\right|^2
=
{1 \over 2^{2(J-L)}} \prod_{j=L}^{J-1} \left|H(2^j\omega)\right|
\cdot \left(\sum_{k=0}^{2^{J-L}-1}\prod_{j=L}^{J-1}
\left|\widetilde{H}(2^j(\omega+k\pi/2^{J-L-1}))\right|\right).
\label{eq:DFTPL}
\end{equation}
In this expression $H(\omega)$ represents the Fourier transform of one row of
the one step matrix $H_j$. In the equidistant settings, all $H_j$
coincide with a single Toeplitz matrix, for which one row suffices to
characterise the complete multiscale transform. Obviously, the same definition
holds for $\widetilde{H}(\omega)$.

The Fourier analysis cannot be applied to non-equidistant observations, and so
it cannot be used to find the best $\widetilde{H}_j$ for a given sequence of
$H_j$.
Still under the assumption that the covariance matrix of $\vec{\varepsilon}$ is
$\Sigma_{\vec{\varepsilon}} = \sigma^2 I$, we find that the covariance matrix
after projection equals
$\Sigma_{\mathrm{P}_L\vec{\varepsilon}} = \sigma^2\mathrm{P}_L\mathrm{P}_L^T$.
Since $\mathrm{P}_L\mathrm{P}_{L\perp} = \mathrm{P}_{L\perp}$, and
$\mathrm{P}_{L\perp}$ is symmetric and idempotent, it holds that 
$\mathrm{P}_L\mathrm{P}_L^T-\mathrm{P}_{L\perp}\mathrm{P}_{L\perp}^T =
\mathrm{P}_L(I-\mathrm{P}_{L\perp})\mathrm{P}_L^T$. The matrix
$I-\mathrm{P}_{L\perp}$ is positive semi-definite, so all diagonal elements of
$\mathrm{P}_L\mathrm{P}_L^T-\mathrm{P}_{L\perp}\mathrm{P}_{L\perp}^T$ must be
non-negative, which reads as $\var(\mathrm{P}_L\vec{\varepsilon}) \geq
\var(\mathrm{P}_{L\perp}\vec{\varepsilon})$, this vector inequality holding
componentwise. This conclusion is confirmed by the observation in Figure
\ref{fig:expvar}. As a result, the variance propagation of a wavelet
decomposition applied to uncorrelated homoscedastic observations is
described by the nonzero eigenvalues of $\mathrm{P}_L\mathrm{P}_L^T$, i.e.,
the nonzero singular values of $\mathrm{P}_L$. As a summary, define the
multiscale variance propagation as follows
\begin{definition} (multiscale variance propagation)
Given a wavelet transform defined by the sequences of forward and inverse 
refinement matrices $H_j$ and $\widetilde{H}_j$, where $j = L,\ldots,J-1$.
Let $W_L$ and $\widetilde{W}_L$ be defined as in (\ref{eq:WL}) and
(\ref{eq:WtL}), and let $\mathrm{P}_L = W_L\widetilde{W}_L$.
Then the multiscale variance propagation equals
\begin{equation}
\kappa_F(W_L,\widetilde{W}_L)
= {\norm{\vec{\sigma}(\mathrm{P}_L)}_2 / \sqrt{n_L}}
= \norm{\mathrm{P}_L}_F / \sqrt{n_L},
\label{def:kappaF}
\end{equation}
where $\vec{\sigma}(\mathrm{P}_L)$ is the vector of singular values of
$\mathrm{P}_L$, while
$n_L$ is the number of columns in $W_L$, which is also the number of nonzero
singular values of $\mathrm{P}_L$. The notation $\norm{\mathrm{P}_L}_F$ stands
for the Frobenius norm of $\mathrm{P}_L$, 
\begin{revision}
whereas 
$\norm{\vec{\sigma}(\mathrm{P}_L)}_2 =
\sqrt{\vec{\sigma}(\mathrm{P}_L)^T\vec{\sigma}(\mathrm{P}_L)}$
is the classical Euclidean vector norm.
\end{revision}
\label{def:multiscalevariancepropagation}
\end{definition}
As an alternative for (\ref{def:kappaF}), the multiscale variance propagation
can also be defined as
\begin{equation}
\kappa_2(W_L,\widetilde{W}_L) = \max(\vec{\sigma}(\mathrm{P}_L))
= \norm{\mathrm{P}_L}_2,
\label{def:kappa2}
\end{equation}
\begin{revision}
where $\norm{\mathrm{P}_L}_2$ denotes the matrix norm induced from the
Euclidean vector norm, which is equal to the largest singular value of the
matrix.
\end{revision}

Using the perfect reconstruction property that $\widetilde{W}_LW_L = I_{n_L}$,
it can be proven that all singular values of
$\mathrm{P}_L = W_L\widetilde{W}_L$ must be
either zero or greater than one. As a consequence, it holds that
$\kappa_F(W_L,\widetilde{W}_L) \geq 1$, and $\kappa_F(W_L,\widetilde{W}_L) = 1$
if and only if $\mathrm{P}_L$ is an orthogonal projection. 

Given the refinement and detail matrices $\widetilde{H}_L^{[0]}$ and
$\widetilde{G}_L$ that result from the factorisation in Section
\ref{subsec:Bspline2scales},
we want to design a sparse update matrix $U_L$ that minimises the Frobenius
norm of $\mathrm{P}_L$, under the constraint that it also preserves a few, say
$p$, primal vanishing moments.
The index $L$ refers here to the coarsest scale up to the current stage in the
design, but $L$ may turn out to be an intermediate scale in the eventual
transform.

We first fix which elements of $U_L$ will be nonzero. The number of non-zeros
must be large enough to cover the $p$ vanishing moments and allow a further
optimisation. For a given element $U_{L;r,s}$ this means that we either choose
it to be zero, or we optimise its value. We therefore compute the derivative of
the Frobenius norm with respect to $U_{L;r,s}$.
Since 
\(
W_L\widetilde{W}_L 
=
W_L(\widetilde{H}_L^{[0]T}+U_L\widetilde{G}_L^T)\widetilde{W}_{L+1},
\)
we find that
\[
{\partial (W_L\widetilde{W}_L)_{i,j} \over \partial U_{L;r,s}}
=
{\partial (W_LU_L\widetilde{G}_L^T\widetilde{W}_{L+1})_{i,j} \over \partial
U_{L;r,s}}
=
W_{L;i,r}(\widetilde{G}_L^T\widetilde{W}_{L+1})_{s,j}.
\]
This we use in
\begin{eqnarray*}
{\partial \norm{W_L\widetilde{W}_L}_F^2 \over \partial U_{L;r,s}}
& = &
{\partial \sum_{i=1}^n\sum_{j=1}^n (W_L\widetilde{W}_L)_{i,j}^2 \over \partial
U_{L;r,s}}
=
\sum_{i=1}^n\sum_{j=1}^n 2(W_L\widetilde{W}_L)_{i,j} 
{\partial (W_L\widetilde{W}_L)_{i,j} \over \partial U_{L;r,s}}
\\
& = &
2\sum_{i=1}^n\sum_{j=1}^n
W_{L;i,r}(W_L\widetilde{W}_L)_{i,j}(\widetilde{G}_L^T\widetilde{W}_{L+1})_{s,j}
=
\left(W_L^TW_L\widetilde{W}_L(\widetilde{G}_L^T\widetilde{W}_{L+1})^T\right)_{r,s}.
\end{eqnarray*}
In this expression, the matrix $\widetilde{W}_L =
(\widetilde{H}_L^{[0]T}+U_L\widetilde{G}_L^T)\widetilde{W}_{L+1}$ depends on
the update matrix $U_L$ that we want to optimise. 

Most applications require at least one primal vanishing moment. Therefore, the
moment equation (\ref{eq:Uvanmom}) is imposed as a constraint in the
optimisation process. Using a vector of Lagrange multipliers for each moment,
the objective function $K_L(U_L,\vec{\lambda}_{L;m})$ can be written as
\[
K_L(U_L,\vec{\lambda}_{L;m}) = \norm{W_L\widetilde{W}_L}_F^2 + \sum_{m=0}^{p-1}
\vec{\lambda}_{L;m}^T (O_{L;m} - U_L^T M_{L;m}).
\]
Given a choice $\mathcal{I}_L = \{(r,s) \in \{1,\ldots,n_L\} \times
\{1,\ldots,n_{L+1}-n_L\}|U_{L;r,s} \neq 0$, the $n_{U_L} = \# \mathcal{I}_L$
equations for the optimisation follow from
\(
{\partial K_L \over \partial U_{L;r,s}} = 0,
\)
\begin{equation}
\left(W_L^TW_LU_L\widetilde{G}_L^T\widetilde{W}_{L+1}\widetilde{W}_{L+1}^T\widetilde{G}_L^T\right)_{r,s}
=
-\left(W_L^TW_L\widetilde{H}_L^{[0]T}\widetilde{W}_{L+1}\widetilde{W}_{L+1}^T\widetilde{G}_L^T\right)_{r,s}
+ \sum_{m=0}^{p-1} \lambda_{L;m;s} M_{L;m;r}.
\label{eq:constrainedFrobeniusoptimisation}
\end{equation}
These equations are completed by $(n_{L+1}-n_L)p$ moment equations
(\ref{eq:Uvanmom}).

The set of non-zeros in $U_L$ is fixed by the user before the optimisation is
carried out. For obvious reasons, its cardinality $n_{U_L}$ should be large
enough to cope with the moment constraints, i.e., $n_{U_L} \geq
(n_{L+1}-n_L)p$. The easiest option is to let $\mathcal{I}_L$ be the collection
of $\nu$ main and side diagonals. The diagonals have lengths equal to
$n_{L+1}-n_L,n_L-1,n_{L+1}-n_L-1,n_L-2,n_{L+1}-n_L-2,\ldots$, summing up to
$n_{U_L} = 
\floor{\nu/2}n_L+\ceil{\nu/2}(n_{L+1}-n_L)-\floor{\nu/2}\ceil{\nu/2}$. Since
$n_{L+1}-n_L$ is the number of odds at scale $j+1$, we have $n_{L+1}-n_L \in
\{n_L-1,n_L\}$, and so 
$n_{U_L} = \nu(n_{L+1}-n_L)-\floor{\nu/2}\ceil{\nu/2}+r_L\floor{\nu/2}$, where
$r_L = 2n_L-n_{L+1} \in \{0,1\}$. It follows that $n_{U_L} \geq (n_{L+1}-n_L)p$
is not possible for $\nu$ equal to $p$, unless $p=1$. In other words, due to
boundary effects, $p$
vanishing moments require more than $p$ diagonals in the final update. 

\begin{figure}[t!]
\includegraphics[width=0.89\textwidth]{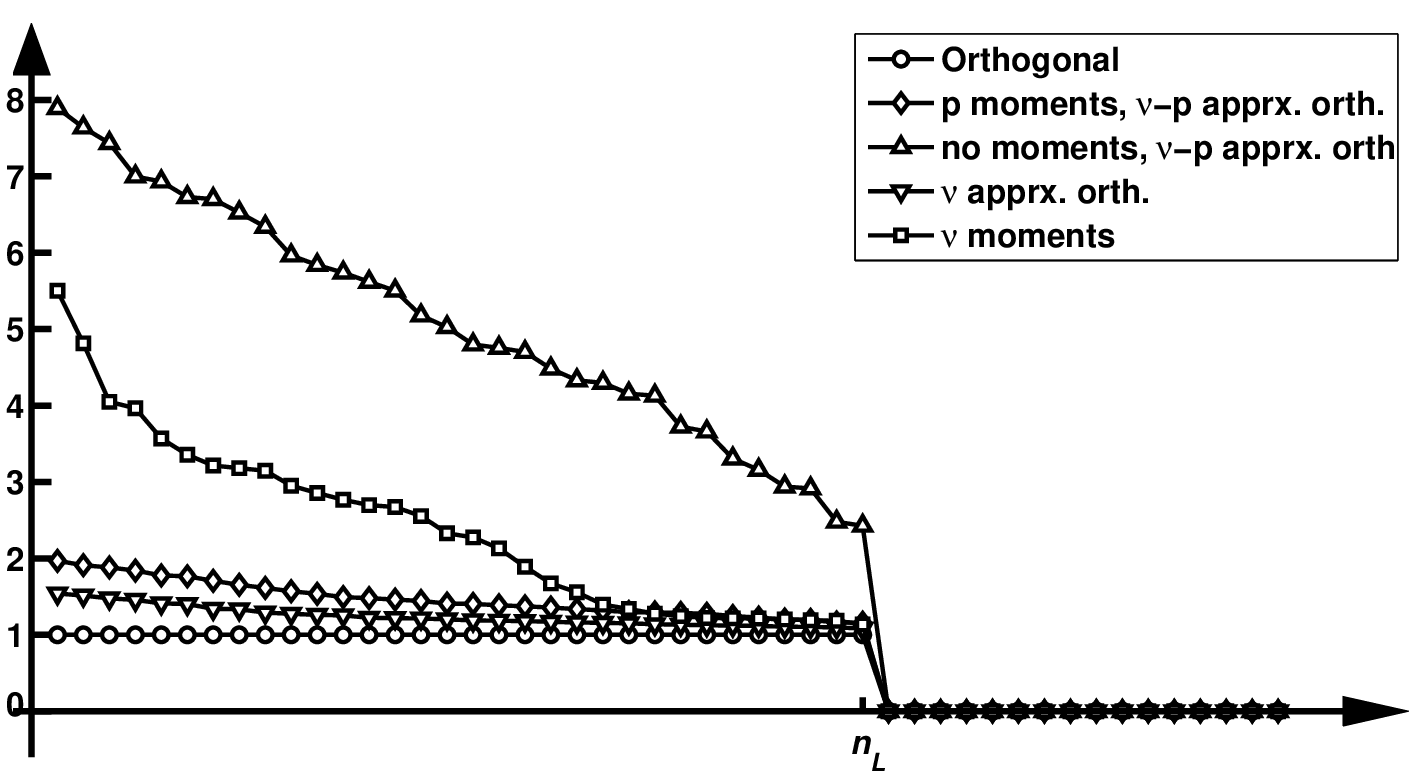}
\caption{Singular values in descending order for projections
$W_L\widetilde{W}_L$ onto a coarse scale B-spline basis, using several forward
transforms $\widetilde{W}_L$. In all examples, $\nu=4$ and $p=2$. The
orthogonal $\widetilde{W}_L$ has no band structure. All other $\widetilde{W}_L$
are quadri-diagonal, except for the upper curve, marked with the $\Delta$ 
signs, which corresponds to a bidiagonal $\widetilde{W}_L$.}
\label{fig:SVDprojBspline}
\end{figure}

As an example, depicted in Figure \ref{fig:SVDprojBspline}, we compute the
singular values of the projections $W_L\widetilde{W}_L$ using cubic B-splines,
for the analysis of data that are defined on $n_J = 1000$ fine scale knots. 
The fine scale knots were drawn independently from each other from a uniform
distribution on $[0,1]$. After $J-L=5$ resolution steps,
the coarse scale has a resolution of $n_L = 32$. Therefore, the matrices $W_L$
and $\widetilde{W}_L^T$ have 1000 rows and 32 columns, from where we know that
$\mathrm{rank}(W_L\widetilde{W}_L) = n_L = 32$. Figure \ref{fig:SVDprojBspline}
plots the first 53 singular values of the projections in descending order. It
is no surprise that the 33rd and beyond are all zero. It is no surprise either
that for the orthogonal projection, i.e., for $\widetilde{W}_L =
(W_L^TW_L)^{-1}W_L^T$, all 32 non-zeros are equal to one. The orthogonal
projection matrix $\widetilde{W}_L$ is not sparse in the strict sense, although
the off-diagonal elements show rapid decay. Therefore, the fully orthogonal
projection can be well approximated by a band-limited matrix, using the system
of equations in (\ref{eq:constrainedFrobeniusoptimisation}). Figure 
\ref{fig:SVDprojBspline} examines the quality of the approximations, showing
that a bidiagonal matrix is probably too sparse: even if all free parameters
are spent in the optimisation of the singular values, the largest value is
close to 8. Three other alternatives for the orthogonal projection use four
diagonals. One option is to spend all free parameters on primal vanishing
moments. The corresponding wavelets $\psi_{L,k}(x)$ all have four vanishing
moments, except those near the boundaries of the interval. The singular values
are, however, not controlled, the maximum being close to 6 in this case. The
situation may deteriorate quickly in other settings. When all free entries in
$U_L$ are spent on the minimisation of the singular values, the fourth order
approximation is ready for use in statistical applications. But even when two
vanishing moments are imposed in combination with a minimisation of the
singular values, these values are kept low enough for use in statistical
applications. Other experiments, not shown here, confirm that imposing one or
two vanishing moments in a scheme that otherwise concentrates on the singular
values, is performing nearly as well in terms of singular values, as a scheme
that has its entire focus on these singular values.

\begin{revision}
As a summary for this section, wavelet transforms for use in statistical
estimation should be as close as possible to being orthogonal, because
reconstructions from non-orthogonal decompositions may suffer from variance
blow up. Orthogonality puts, however, serious limitations on the design of a
wavelet transform. As an example, orthogonal spline wavelets with compact 
support do not exist.
The numerical notions of condition number and Riesz constants in a
biorthogonal transform are not sufficiently adequate for the description of the
variance propagation. The design of the wavelet transform proposed in this
section therefore focusses directly on the variance and it does so by looking
at the combined effect of decomposition and reconstruction.
\end{revision}

\section{Conclusion and outlook}

The first contribution of this paper has been to extend the construction of 
the Cohen-Daubechies-Feauveau B-spline wavelets towards the case of
non-equidistant knots. The new construction is based on the factorisation of
two-scale or refinement equations into lifting steps. An interesting topic for
further research is to investigate the same methodology for other spline
wavelets.

The second contribution has been the design of a numerical method to find
fine scaling coefficients from function values, in a way that performs well
when the function is observed with noise. Future research could improve the
method by making it adaptive to jumps or other singularities. In the
neighbourhood of jumps, the benefit from a perfect reconstruction of power
functions is limited, so a local relaxation of these conditions may yield
sharper reconstructions.

Finally, the third contribution has been a modification of the
Cohen-Daubechies-Feauveau wavelet transform for specific use in statistical
applications, making sure to control the propagation of the variance on the
wavelet coefficients at successive scales. 
\begin{revision}
The focus on the variance propagation in the design of the transform
allows us to relax the stringent orthogonality condition and to construct a 
basis of compactly supported spline wavelets in which non-linear processing
leads to a reconstruction with a well controlled bias-variance balance.
\end{revision}

\subsection*{Acknowledgement}

Research support by the IAP research network grant nr. P7/06 of the Belgian
government (Belgian Science Policy) is gratefully acknowledged.

\appendix
\section*{Appendices}

\section{Available software - reproducible research}
\label{app:software}

All transforms presented in this article have been implemented in the latest
version of \texttt{ThreshLab}, a Matlab\textregistered software package
available for download from \\
\url{http://homepages.ulb.ac.be/~majansen/software/threshlab.html}.

The forward and inverse B-spline wavelet transforms are carried out by the
routines \texttt{FWT\_2Gspline.m} and \texttt{IWT\_2Gspline.m}. 
Several alternatives for the retrieval of appropriate fine scaling coefficients
from noisy observations, explained in Section \ref{sec:apprxerr}, are
implemented in \texttt{finescaleBsplinecoefs.m}. The use of this routine is
illustrated in \texttt{illustratefinescalecoefs.m}.

In particular, the experiment in Figure \ref{fig:expvar}
is set up in the routine \texttt{illustratevarianceprojection.m}.
The singular value plots in Figure \ref{fig:SVDprojBspline}
have been generated using \texttt{illustrate\_updateLSapprxprimmom.m}.

\begin{longversion}
\section{Proofs for Theorems \ref{thm:Bsplinebasis},\ref{thm:x^qcoefBspline},
and \ref{thm:Bsplinesfrompolynomialreproduction}}
\label{app:proof:thms:Bsplinebasis,x^qcoefBspline,Bsplinesfrompolynomialreproduction}

The proofs are given for the sake of self-containment. The results are
well-known in literature.

\subsection{Proof of Theorem \ref{thm:Bsplinebasis}}
\label{app:proof:thm:Bsplinebasis}

The summation in (\ref{eq:Bsplinebasis}) contains
$n_j-\ceil{\widetilde{p}/2}-\floor{\widetilde{p}/2} = n_j - \widetilde{p}$
B-spline functions. All B-splines have mutually
unequal supports, and are thus linearly independent.
On the other hand, $I_j$ consists of $n_j-2\widetilde{p}+1$ subintervals
$[x_{j,k},x_{j,k+1})$, with $\widetilde{p}$ degrees of freedom on each of them.
The continuous derivatives in each interior knot consume
$(n_j-2\widetilde{p})\times(\widetilde{p}-1)$ of these degrees of freedom,
leaving us with a vector space of dimension $n_j - \widetilde{p}$.
\qed

\subsection{Proof of Theorem \ref{thm:x^qcoefBspline}}
\label{app:proof:thm:x^qcoefBspline}

\begin{enumerate}
\setcounter{enumi}{1}
\item
Expression (\ref{eq:polynomialcoefsBspline}) follows from
(\ref{eq:derivativeinBsplinebasis}), applied for $f_j(x) = x^q$. The left hand
side of (\ref{eq:derivativeinBsplinebasis}) can then be written as
\[
qx^{q-1} = q\sum_{k \in \ZZ} 
\widetilde{x}_{j,k-\widetilde{r}'}^{\left[\widetilde{p}-1,q-1\right]}
\varphi_{j,k-\widetilde{r}'}^{\left[\widetilde{p}-1\right]}(x),
\]
whereas the right hand side becomes
\[
qx^{q-1} = (\widetilde{p}-1) \sum_{k\in \ZZ} 
         {\widetilde{x}_{j,k}^{\left[\widetilde{p},q\right]} -
          \widetilde{x}_{j,k-1}^{\left[\widetilde{p},q\right]} \over
           x_{j,k+\ceil{\widetilde{p}/2}-1}-x_{j,k-\floor{\widetilde{p}/2}}}
           \varphi_{j,k-\widetilde{r}'}^{\left[\widetilde{p}-1\right]}(x).
\]
Identification of the terms in the decomposition leads to
(\ref{eq:polynomialcoefsBspline}).
\item
Next, it can be verified that all solutions for
$\widetilde{x}_{j,k}^{\left[\widetilde{p},1\right]}$ in
(\ref{eq:polynomialcoefsBspline}) must satisfy
\[
\widetilde{x}_{j,k}^{\left[\widetilde{p},1\right]}
=
{1 \over \widetilde{p}-1}
\sum_{i=1-\floor{\widetilde{p}/2}}^{\ceil{\widetilde{p}/2}-1} x_{j,k+i}
-
{1 \over \widetilde{p}-1}
\sum_{i=1-\floor{\widetilde{p}/2}}^{\ceil{\widetilde{p}/2}-1} x_{j,i} +
s_{j,0}^{\left[\widetilde{p},1\right]}.
\]
On the other hand, $\widetilde{x}_{j,k}^{\left[\widetilde{p},1\right]}$ must be
independent from the $\widetilde{p}-1$ knots $x_{j,i}$ around $x_{j,0}$
if $k > \widetilde{p}$. So take $x_{j,i}$ symmetric around $x_{j,0}=0$.
Then, obviously,
\[
{1 \over \widetilde{p}-1}
\sum_{i=1-\floor{\widetilde{p}/2}}^{\ceil{\widetilde{p}/2}-1} x_{j,i} = 0.
\]
On the other hand, $s_{j,0}^{\left[\widetilde{p},1\right]} = 0$, as the
corresponding basis function is even. This leads to
(\ref{eq:linearcoefsBspline}).
\item
The proof for (\ref{eq:x^(p-1)coefsBspline}) is similar to that for
(\ref{eq:linearcoefsBspline}), following an induction argument on
$\widetilde{p}$.
\end{enumerate}
\qed

\subsection{Proof of Theorem \ref{thm:Bsplinesfrompolynomialreproduction}}
\label{app:proof:thm:Bsplinesfrompolynomialreproduction}

Consider $x \in (x_{j,k},x_{j,k+1}) \subset I_j$, then
\(
x \in S_{j,m} \Leftrightarrow l_m < k+1 \mbox{ and } r_m > k
\Leftrightarrow
k-\ceil{\widetilde{p}/2}+1 \leq m \leq k+\floor{\widetilde{p}/2}.
\)
Then in $x$, it holds that
\[
\left[1 \,\, x \,\, x^2 \,\, \ldots x^{\widetilde{p}-1}\right] =
\left[\varphi_{j,k-\ceil{\widetilde{p}/2}+1}^{\left[\widetilde{p}\right]}(x)
\,\, \ldots\,\,
\varphi_{j,k+\floor{\widetilde{p}/2}}^{\left[\widetilde{p}\right]}(x)\right]
\,\widetilde{X},
\]
where
\(
\widetilde{X}_{m-k+\ceil{\widetilde{p}/2},q+1} =
\widetilde{x}_{j,m}^{\left[\widetilde{p},q\right]}.
\)
It can be verified that $\widetilde{X}$ is non-singular, so we find an
expression for the basis functions on the subinterval $(x_{j,k},x_{j,k+1})$,
\begin{equation}
\left[\varphi_{j,k-\ceil{\widetilde{p}/2}+1}^{\left[\widetilde{p}\right]}(x)
\,\, \ldots \,\,
\varphi_{j,k+\floor{\widetilde{p}/2}}^{\left[\widetilde{p}\right]}(x)\right]
=
\left[1 \,\, x \,\, x^2 \,\, \ldots x^{\widetilde{p}-1}\right]
\,\widetilde{X}^{-1},
\label{eq:phifromxtilde}
\end{equation}
all other basis functions being zero on this subinterval. It follows that on
each subinterval, the basis functions must be polynomials, and the polynomials
are uniquely defined by (\ref{eq:phifromxtilde}). As from Theorem
\ref{thm:x^qcoefBspline}, we know that the power functions have the
coefficients $\widetilde{x}_{j,m}^{\left[\widetilde{p},q\right]}$ in a B-spline
basis, the polynomials defined by (\ref{eq:phifromxtilde}) must coincide with
the B-splines on that interval.
\qed

\end{longversion}

\section{Constructive proof for the factorisation in Theorem
\ref{thm:factorisation}}
\label{app:proof:thm:factorisation}

The factorisation is based on the approach in
\citet{daubechies98:factorlifting}. Because of the non-equidistant knots, it
proceeds column by column in the refinement matrix. The proof also makes use of
the band structure of a matrix. For the sake of simplicity, and without any
impact on the applicability of the argument, we ignore occasional zeros within
the nonzero band of $H_{j,e}^{[s]}$ or $H_{j,o}^{[s]}$. In each column of 
these matrices, we take the same number of rows into consideration, equal to
the bandwidth of the whole matrix, even if that particular row has actually
less non-zeros.

Consider first the case where $H_{j,e}^{[s]}$ has a larger bandwidth
than $H_{j,o}^{[s]}$. This is only possible if the first and last nonzero in
each column of $H_j$ is situated on an even row. Denote by $2k_1(\ell)$ the
first nonzero row in column $\ell$ of $H_j$ and by $2k_2(\ell)$ the last
nonzero row in the same column.
Then the $\ell$th column of (\ref{eq:primallifting}) reads as
\begin{equation}
H_{j,2k,\ell}^{[s]} = H_{j,2k,\ell}^{[s+1]} - \sum_{m = k_1(\ell)}^{k_2(\ell)}
U_{j,k,m}^{[s+1]} H_{j,2m+1,\ell}^{[s]},
\label{eq:factorU1}
\end{equation}
We assume that $k_1(\ell)$ and $k_2(\ell)$ are strictly increasing
functions, otherwise a more careful design of the update step is needed.
Let $\ell_2(k)$ be the inverse of $k_1(\ell)$, i.e.,
$\ell_2(k_1(\ell)) = \ell$. In words, $\ell_2(k)$ is the last column on row 
$k$ with a nonzero element.
In a similar way, let $\ell_1(k)$ be the inverse of $k_2(\ell)$.
In a given column $\ell$, we impose for
$k=k_1(\ell)$ and for $k=k_2(\ell)$ that $H_{j,2k,\ell}^{[s+1]} = 0$, so that
the number of non-zeros in column $\ell$ of $H_{j,e}^{[s+1]}$ equals
$k_2(\ell)-k_1(\ell)-1$.
Note that the number of nonzero rows equals $k_2(\ell)-k_1(\ell)+1$ for
$H_{j,e}^{[s]}$ and $k_2(\ell)-k_1(\ell)$ for $H_{j,o}^{[s]}$.
For $k=k_1(\ell)$, and for $k=k_2(\ell)$, we obtain
the two equations of the form
\begin{equation}
H_{j,2k,\ell}^{[s]} = - \sum_{m=k_1(\ell)}^{k_2(\ell)-1}U_{j,k,m}^{[s+1]} 
H_{j,2m+1,\ell}^{[s]}.
\label{eq:factorU2}
\end{equation}
These two equations contain as unknowns two partial rows in matrix
$U_j^{[s+1]}$.
Instead of solving the two equations for fixed $\ell$, we consider the two
equations for given $k$, by looking at $\ell = \ell_2(k)$ and
$\ell = \ell_1(k)$.

For $k=k_1(\ell)$, i.e., $\ell = \ell_2(k)$, we set $U_{j,k,m}^{[s+1]} = 0$
for all $m = k+1,\ldots,k_2(\ell_2(k))-1$, while for $m=k$, we must then take
\begin{equation}
U_{j,k,k}^{[s+1]} = -H_{j,2k,\ell_2(k)}^{[s]}/H_{j,2k+1,\ell_2(k)}^{[s]},
\label{eq:designUdiag}
\end{equation}
in order to satisfy (\ref{eq:factorU2}) for $k=k_1(\ell)$ or $\ell = \ell_2(k)$.

For $k=k_2(\ell)$, i.e., $\ell = \ell_1(k)$, we set $U_{j,k,m}^{[s+1]} = 0$
for all $m = k_1(\ell_1(k)),\ldots,k-2$, while for $m=k-1$, we must then take
\begin{equation}
U_{j,k,k-1}^{[s+1]} = -H_{j,2k,\ell_1(k)}^{[s]}/H_{j,2k-1,\ell_1(k)}^{[s]},
\label{eq:designUlowerdiag}
\end{equation}
in order to satisfy (\ref{eq:factorU2}) for $k=k_2(\ell)$ or $\ell = \ell_1(k)$.

Once the diagonal and the lower diagonal of $U_j^{[s+1]}$ has been found, all
the other entries of this matrix can be filled with zeros. This leaves us with
$k_2(\ell)-k_1(\ell)-1$ equations of the form (\ref{eq:factorU1}). Each of
these equations allows us to find exactly one element in column $\ell$ of
$H_{j,e}^{[s+1]}$.

The case where the bandwidth of $H_{j,o}^{[s]}$ is larger than that of
$H_{j,e}^{[s]}$ is treated in a similar way.
The case where the bandwidths of $H_{j,e}^{[s]}$  and $H_{j,o}^{[s]}$ are equal
can be reduced to the first case if we artificially increase the bandwidth of
$H_{j,e}^{[s]}$, by taking an additional zero into account in each column of
$H_{j,e}^{[s]}$.
\qed


\section{Proof of Proposition \ref{prop:liftedcoefficient}}
\label{app:proof:prop:liftedcoefficient}

First, it can be checked that (\ref{eq:liftingoutandin}) holds for any
$\widetilde{p}$ and for $q=1$ and for $q=\widetilde{p}-1$. More precisely, a
bit of calculations show that from (\ref{eq:linearcoefsBspline}), it follows
indeed that,
\begin{equation}
\widetilde{t}_{j,k}^{\left[\widetilde{p},1,L1\right]}
=
{t_{j,k+\widetilde{p}-1}-t_{j,k+1} \over t_{j,k+\widetilde{p}}-t_{j,k}}
\,
\left(
t_{j,k} + \sum_{i=2}^{\widetilde{p}-2} t_{j,k+i} + t_{j,k+\widetilde{p}}
\right).
\end{equation}
In a similar way, starting from (\ref{eq:x^(p-1)coefsBspline}), we arrive at
\begin{equation}
\widetilde{t}_{j,k}^{\left[\widetilde{p},\widetilde{p}-1,L1\right]}
=
{t_{j,k+\widetilde{p}-1}-t_{j,k+1} \over t_{j,k+\widetilde{p}}-t_{j,k}}
\,
\left(
t_{j,k} \cdot \prod_{i=2}^{\widetilde{p}-2} t_{j,k+i} \cdot
t_{j,k+\widetilde{p}}
\right).
\end{equation}

As a result, for $\widetilde{p}=1,2$ the result holds for all
$q=0,\ldots,\widetilde{p}-1$. This is the basis for the subsequent induction
argument.

Using the recursion in (\ref{eq:polynomialcoefsBsplineN}), we can find the
coefficient $\widetilde{t}_{j,k}^{\left[\widetilde{p},q,L1b\right]}$ in two
steps. First we define $\widetilde{t}_{j,k}^{\left[\widetilde{p},q,L1a\right]}$
as the coefficient that results from taking out $t_{j,k+1}$ from
$\widetilde{t}_{j,k}^{\left[\widetilde{p},q\right]}$ and replacing it by
$t_{j,k}$. This intermediate coefficient equals
\[
\widetilde{t}_{j,k}^{\left[\widetilde{p},q,L1a\right]} = 
\widetilde{t}_{j,k}^{\left[\widetilde{p},q\right]} +
{q \over \widetilde{p}-1}\,
\widetilde{t}_{j,k+1}^{\left[\widetilde{p}-1,q-1\right]}\,
\left(t_{j,k}-t_{j,k+1}\right).
\]
We define a similar coefficient for $\widetilde{p}-1$ and $q-1$.
\[
\widetilde{t}_{j,k}^{\left[\widetilde{p}-1,q-1,L1a\right]} = 
\widetilde{t}_{j,k}^{\left[\widetilde{p}-1,q-1\right]} +
{q-1 \over \widetilde{p}-2}\,
\widetilde{t}_{j,k+1}^{\left[\widetilde{p}-2,q-2\right]}\,
\left(t_{j,k}-t_{j,k+1}\right).
\]
Next, we take out $t_{j,k+\widetilde{p}-1}$ from
$\widetilde{t}_{j,k}^{\left[\widetilde{p},q,L1a\right]}$ and replace it with
$t_{j,k+\widetilde{p}}$. In order to apply (\ref{eq:polynomialcoefsBsplineN}),
we need the power coefficient of degree
$\widetilde{p}-1$ for $q-1$ based on all remaining knots in
$\widetilde{t}_{j,k}^{\left[\widetilde{p},q,L1a\right]}$. This is exactly
$\widetilde{t}_{j,k}^{\left[\widetilde{p}-1,q-1,L1a\right]}$. We thus find
\[
\widetilde{t}_{j,k}^{\left[\widetilde{p},q,L1b\right]} =
\widetilde{t}_{j,k}^{\left[\widetilde{p},q,L1a\right]} +
{q \over \widetilde{p}-1}\,
\widetilde{t}_{j,k}^{\left[\widetilde{p}-1,q-1,L1a\right]}
\left(t_{j,k+\widetilde{p}}-t_{j,k+\widetilde{p}-1}\right).
\]
Using the expressions above for
$\widetilde{t}_{j,k}^{\left[\widetilde{p},q,L1a\right]}$ and for
$\widetilde{t}_{j,k}^{\left[\widetilde{p}-1,q-1,L1a\right]}$, 
this can further be developed into
\begin{equation}
\begin{array}{rl}
\displaystyle
\widetilde{t}_{j,k}^{\left[\widetilde{p},q,L1b\right]}
=
\widetilde{t}_{j,k}^{\left[\widetilde{p},q\right]} +
{q \over \widetilde{p}-1}
&
\!\!\!\!
\left[
\widetilde{t}_{j,k}^{\left[\widetilde{p}-1,q-1\right]}
(t_{j,k}+t_{j,k+\widetilde{p}}-t_{j,k+1}-t_{j,k+\widetilde{p}-1})
\right.
\\
&
\displaystyle
\left.
+
{q-1 \over \widetilde{p}-2}\,
\widetilde{t}_{j,k+1}^{\left[\widetilde{p}-2,q-2\right]}
(t_{j,k}-t_{j,k+1})(t_{j,k+\widetilde{p}}-t_{j,k+1})
\right].
\end{array}
\label{eq:2out2inbyrecursion}
\end{equation}

Expression (\ref{eq:2out2inbyrecursion}) can be substituted in the right hand
side of (\ref{eq:liftingoutandin}). We now work on
$\widetilde{t}_{j,k}^{\left[\widetilde{p},q,L1\right]}$ in the left hand side,
starting from its definition in (\ref{eq:liftedpowercoef}), and using once more 
the recursion (\ref{eq:polynomialcoefsBsplineN}).
We find
\begin{eqnarray*}
\widetilde{t}_{j,k}^{\left[\widetilde{p},q,L1\right]}
& = &
{t_{j,k+\widetilde{p}-1}-t_{j,k+1} \over t_{j,k+\widetilde{p}}-t_{j,k}}
\,\,\widetilde{t}_{j,k}^{\left[\widetilde{p},q\right]} 
- {q \over \widetilde{p}-1} \,\, {1 \over t_{j,k+\widetilde{p}}-t_{j,k}}\cdot
\\
&&
\left[
\widetilde{t}_{j,k+1}^{\left[\widetilde{p}-1,q-1\right]}
(t_{j,k+\widetilde{p}}-t_{j,k+1})(t_{j,k+1}-t_{j,k})
+
\widetilde{t}_{j,k}^{\left[\widetilde{p}-1,q-1\right]}
(t_{j,k}-t_{j,k+\widetilde{p}-1})
(t_{j,k+\widetilde{p}}-t_{j,k+\widetilde{p}-1})
\right].
\end{eqnarray*}
For the last factor in this expression, we use again the recursion
(\ref{eq:polynomialcoefsBsplineN}), this time for
\[
\widetilde{t}_{j,k+1}^{\left[\widetilde{p}-1,q-1\right]}
=
\widetilde{t}_{j,k}^{\left[\widetilde{p}-1,q-1\right]}
+
{q-1 \over \widetilde{p}-2}\,
\widetilde{t}_{j,k+1}^{\left[\widetilde{p}-2,q-2\right]}\,
\left(t_{j,k+\widetilde{p}-1}-t_{j,k+1}\right).
\]
Substitution of this recursion, followed by straightforward algebraic
manipulation, amounts to (\ref{eq:liftingoutandin}), thereby completing the
proof.

\bibliographystyle{plainnat}

\subsubsection*{Correspondence address}

\begin{singlespace}
Maarten Jansen\\
Departments of Mathematics and Computer Science,\\
Universit\'e Libre de Bruxelles\\
Boulevard du Triomphe\\
Campus Plaine, CP213\\
B-1050 Brussels - Belgium \\
\texttt{maarten.jansen@ulb.ac.be}
\end{singlespace}

\end{document}